%
%
%
%
%
%
%
%
%
%
\scrollmode
\magnification=\magstep1
\parskip=\smallskipamount

\def\demo#1:{\par\medskip\noindent\it{#1}. \rm}
\def\ni{\noindent}               
\def\ll{\leftline}
\def\cl{\centerline}

\def\begin{\ll{}\vskip 10mm \nopagenumbers}  
\def\pn{\footline={\hss\tenrm\folio\hss}}   
\def\ii#1{\itemitem{#1}}

%
%
\outer\def\beginsection#1\par{\bigskip
  \message{#1}\leftline{\bf\&#1}
  \nobreak\smallskip\vskip-\parskip\noindent}

%
%
\outer\def\proclaim#1:#2\par{\medbreak\vskip-\parskip
    \noindent{\bf#1.\enspace}{\sl#2}
  \ifdim\lastskip<\medskipamount \removelastskip\penalty55\medskip\fi}

\def\endpr{\hfill $\spadesuit$ \medskip}

%
%
%
%


%
%
%
%

\def\C{{\bf C}}

\def\R{{\bf R}}

%
%
%
%

\def\cC{{\cal C}}

%
%
%
\def\a{\alpha}
\def\b{\beta}

\def\d{\delta}
\def\e{\epsilon}

\def\l{\lambda}

\def\c{\chi}


%
%
%
%
\def\bar{\overline}              
\def\bs{\backslash}              
\def\di{\partial}                
\def\dibar{\bar\partial}         
\def\hra{\hookrightarrow}

%
%
\def\cn{{\bf C}^n}

\def\c*{{\bf C}^*}

%
%
\def\holo{holomorphic}                   
\def\nbd{neighborhood}                   
\def\psc{pseudoconvex}                   
\def\spsc{strongly\ pseudoconvex}        
\def\ra{real-analytic}                   
\def\psh{plurisubharmonic}               
\def\spsh{strongly\ plurisubharmonic}
\def\tr{totally real}                    
\def\hc{holomorphically convex}          
\def\ss{\subset\!\subset}                

\def\supp{{\rm supp}\,}                  

\def\hvf{holomorphic vector field}
\def\hvb{holomorphic vector bundle}
\def\vf{vector field}

\def\pb{\psc\ bump}

\def\Aut{{\rm Aut}}                         

\begin
\cl{\bf OKA'S PRINCIPLE FOR HOLOMORPHIC FIBER BUNDLES}
\medskip
\cl{\bf WITH SPRAYS}
\bigskip\medskip
\cl{Franc Forstneri\v c and Jasna Prezelj}
\bigskip\medskip\rm

\beginsection 1. The Oka-Grauert principle.

The subject of this paper is the homotopy principle, also called
the h-principle or the Oka-Grauert principle, concerning sections of
certain \holo\ fiber bundles on Stein manifolds. We give a proof of a theorem
of Gromov (1989) from sec.\ 2.9 in [Gro]; see theorems 1.3 and 1.4
below. This result, which extends the work of H.\ Grauert from 1957
([Gr3], [Gr4], [Car]), has been used in the proofs of the embedding
theorem for Stein manifolds into Euclidean spaces of minimal dimension
[EGr], [Sch].

\proclaim 1.1 Definition:
Let $h\colon Z\to X$ be a holomorphic mapping of complex manifolds.
A {\it section} of $h$ is any map $f\colon X\to Z$ such that
$h\circ f$ is the identity on $X$. We say that sections of $h$ satisfy
the {\bf h-principle} (or the {\bf Oka-Grauert principle}) if each
continuous section $f_0 \colon X\to Z$ can be deformed to a
holomorphic section $f_1\colon X\to Z$ through
a continuous one parameter family (a homotopy) of continuous sections
$f_t\colon X\to Z$ ($0\le t\le 1$), and any two holomorphic sections
$f_0,f_1\colon X\to Z$ which are homotopic through continuous sections
are also homotopic through holomorphic sections. If this holds
for a trivial bundle $Z=X\times F\to X$, we say that maps $X\to F$
satisfy the h-principle.

%
%
\proclaim 1.2 Definition: {\rm (Gromov [Gro])}
A (dominating) {\bf spray} on a complex manifold $F$ is a \hvb\
$p\colon E\to F$, together with a holomorphic map $s\colon E\to F$, such that
$s$ is the identity on the zero section $F\subset E$, and for each $x\in F$
the derivative $Ds(x)$ maps $E_x$ (which is a linear subspace of $T_x E$)
surjectively onto $T_x F$.

The following result can be found in sec.\ 2.9 of [Gro].

%
%
\proclaim 1.3 Theorem: 
If $F$ is a complex manifold which admits a spray, then the sections
of any locally trivial \holo\ fiber bundle with fiber $F$ over any Stein
manifold satisfy the h-principle.  In particular, mappings from Stein
manifolds into $F$ satisfy the h-principle.

Stronger results are given in theorem 1.4 and corollary 1.5 below.
In the sequel [FP] to this paper we give a proof of Gromov's
Main Theorem ([Gro], sect.\ 4.5) to the effect that the h-principle
holds for sections of \holo\ submersions $h\colon Z\to X$, where
$X$ is Stein and each point $x\in X$ has a \nbd\ $U \subset X$
such that $Z|U=h^{-1}(U)$ admits a fiber-spray (see def.\ 3.1 below).

For non-specialists we recall that a complex manifold is called
{\it Stein} (after Karl Stein, 1951 [Ste]) if it has `plenty' of
global \holo\ functions. For the precise definition and properties
we refer the reader to the monographs [GRo], [H\"o2] and [GRe].
The most commonly used characterizations are the following.
A complex manifold $X$ is Stein if and only if any of the
following two conditions holds:

\item{--} $X$ can be embedded as a closed complex submanifold of
some complex Euclidean space (the embedding theorem of Remmert,
Bishop, and Narasimhan [GRo, p.\ 224]);

\item{--}  $X$ admits a smooth \spsh\ exhaustion function (Grauert [Gr1]).

\medskip
By a {\it locally trivial \holo\ fiber bundle} with fiber $F$ over a complex
manifold $X$ we mean a bundle obtained by patching the trivial bundles
$U_\a \times F$ over an open covering $\{U_\a \}$ of $X$ by transition
functions of the form
$$
    \phi_{\a,\b} (x,\xi) = \bigl(x, \psi_{\a,\b}(x,\xi)\bigr)
    \qquad (x\in U_{\a,\b}= U_\a \cap U_\b,\ \xi\in F),
$$
where $\psi_{\a,\b} \colon U_{\a,\b} \times F\to F$ is holomorphic
and $\psi_{\a,\b}(x,\cdotp) \in \Aut F$ is an automorpism of $F$ for each
fixed $x$. Of course these transition function must satisfy the
usual compatibility conditions, see for instance [Car].

\pn

Examples of bundles satisfying theorem 1.3 include those
whose fiber is a complex Lie group $L$ or an $L$-homogeneous
space; see examples (A) and (B) below for the existence of spray.
The simplest bundles of this type are the principal holomorphic
$L$-bundles in which the transition maps are given left
multiplications by holomorphic $L$-valued functions
$U_{\a,\b}\to L$, and $L$ acts on each fiber $E_x$ by right
multiplication.  Another example are bundles in which
the transition maps $\psi_{\a,\b}(x,\cdotp)$ are Lie group
automorphisms of $L$. In these cases the h-principle had been
proved by Grauert [Gr3], [Car].

In general the automorphisms $\psi_{\a,\b}(x,\cdotp)$
need not preserve any additional structure the fiber might
have and hence we get much more general fiber bundles.
In this connection we recall an important example of
Demailly [Dem]:

%
%

%
%
\demo Example 1: There exists a locally trivial \holo\ fiber
bundle $h\colon Z \to X$, with base $X$ either $\C$ or the disc
and with fiber $\C^2$,
%
%
such that $Z$ has no holomorphic functions other than those of
the form $g\circ h$, where $g$ is \holo\ on the base $X$.
So $Z$ is not Stein. In particular, $Z$ admits no \hvb\ structure
which gives a negative answer to the question of Gromov [Gro, 2.5.B].
The transition functions used in the construction of $Z$ are (nonlinear)
\holo\ automorphisms of the fiber $\C^2$. Nevertheless the sections of
$Z\to X$ satisfy the h-principle according to theorem 1.3.
\endpr

%
%
%
%

%
%
One of the most important sources of spaces with sprays
is the following [Gro]:
%
%

\medskip \ni \it
\item{(*)} If $F$ is a complex manifold which admits finitely many
$\C$-complete \hvf s $V_j$, $1\le j\le J$, such that the vectors
$V_j(x)$ span the tangent space $T_x F$ at each point $x\in F$,
then $F$ admits a spray. \rm

\medskip
Recall that a \hvf\ is $\C$-complete if for any $x\in F$ the
flow $t\to \phi^t(x)$ of $V$, with $\phi^0(x)=x$,
is defined for all complex values
$t\in \C$ of the time parameter. (For a discussion of
$\R$-completeness versus $\C$-completeness see [Fo1] and [AFR].)
Indeed, if we denote by $\phi_j^t$ the flow of $V_j$ on $F$
and let $s\colon F\times \C^J \to F$ be given by
$$
      s(x;t_1,\ldots,t_J) =
      \phi_1^{t_1}\circ \phi_2^{t_2}\circ\cdots\circ
      \phi_J^{t_J}(x),
$$
we have $s(x;0,\ldots,0)=x$ and
${\di\over \di t_j} s(x;0,\ldots,0)=V_j(x)$ ($x\in F$, $1\le j\le J$).
If these vectors span $T_x F$ for each $x\in F$,
$s$ is a spray on $F$.

\medskip
%
%
We list some examples of complex manifolds with sprays:
%
%

\item{(A)} Any complex Lie group $L$. (Take any finite set of left
or right invariant vector fields on $L$ which span $T_e L=$the
Lie algebra of $L$.)

\item{(B)} Let $L$ be a complex Lie group with Lie algebra
$\ell$. If $L$ acts holomorphically and transitively
on a complex manifold $F$, the map $s\colon F\times \ell \to F$,
$s(x,l)=\exp(l)x$ ($x\in F$, $l\in\ell$), is a spray on $F$.
$F$ is then biholomorphic to the $L$-homogeneous space
$L/H=\{lH\colon l\in L\}$ where $H$ is a complex Lie subgroup
of $L$ (the isotropy group of a point in $Y$).

\item{(C)} $F= \cn \bs \Sigma$, where $\Sigma$ is an algebraic subvariety
of complex codimension at least two (see below).

\medskip
We emphasize that, at the time of this writing,
{\it the known proofs of the embedding theorem for Stein manifolds
(and Stein spaces) to Euclidean spaces of minimal dimension
depend on the validity of the h-principle in the cases (B) and (C)
\rm (see Eliashberg and Gromov [EGr] and Sch\"urmann [Sch]).
\medskip

In case (C) we think of $\Sigma$ as a set which should be avoided
by the image of a map $X\to \cn$. Theorem 1.3 asserts that we can
avoid $\Sigma$ by a holomorphic map $f\colon X\to \cn$ if we can
do so by a continuous map. In this case we can obtain complete
vector fields on $F$ of the form $V(z)=f(\pi(z))v$ ({\it shear fields}),
where $v\in \C^n\bs\{0\}$,
$\pi\colon\cn\to\C^{n-1}$ is a linear projection with $\pi(v)=0$,
and $f\colon \C^{n-1}\to\C$ is an entire function which vanishes
on the projected set $\Sigma'=\pi(\Sigma) \subset \C^{n-1}$. The
flow of $V$, given by $\phi^t(z)=z+tf(\pi(z))v$, fixes $\Sigma$
and hence induces a complete flow on $\C^n\bs \Sigma$.
The projection $\pi$ must be chosen such that it is proper when
restricted to $\Sigma$ to insure that $\pi(\Sigma)=\Sigma'$ is a
proper closed subvariety of $\C^{n-1}$ of codimension at least one.
This is true for most projections, see e.g.\ [Chi].

The argument above works as long as there exist sufficiently
many linear projection $\pi\colon \cn\to\C^{n-1}$ which are proper
when restricted to $\Sigma$. The following example shows that
one cannot completely dispose of the last condition.

%
%
\demo Example 2: For any integer $N>0$ there exist discrete sets
$\Sigma\subset \C^N$  for which there exist no non-degenerate \holo\ maps
$\C^N\to \C^N\bs\Sigma$ [RRu]. Furthermore, for any $1\le n < N$ there
exist proper holomorphic embeddings $\sigma\colon \C^n \to \C^N$ such that
$\C^N\bs \sigma(\C^n)$ admits no non-degenerate holomorphic images
of $\C^{N-n}$ [BFo], [Fo2], and hence
$\C^N\bs \sigma(\C^n)$ admits no spray.
In [FP] we show that {\it the h-principle fails in these cases.}
\medskip

In these constructions,  as in the work of Grauert [Gr3] and
Gromov [Gro], one must consider not a single section at a time
but families of sections, depending continuously on a parameter
$y$ in a compact Hausdorff space $Y$. The basic objects will be
continuous maps $f\colon X\times  Y\to Z$ such that
$f(\cdotp,y)\colon X\to Z$ is a section of $h\colon Z\to X$ for
each fixed $y\in Y$. A {\it homotopy} of such maps is a continuous
map $H\colon X\times Y\times [0,1]\to Z$ such that
$H_t(\cdotp,y)=H(\cdotp,y,t) \colon X\to Z$ is a section of
$h\colon Z\to X$ for all $y\in Y$ and $t\in [0,1]$.

Recall that a compact subset $K \subset X$ is {\it \hc} in $X$ if for
each $x\in X\bs K$ there is a \holo\ function $f$ on $X$ such that
$|f(x)|>\sup_K |f|$. If $X$ is Stein then by the Oka-Weil theorem
each function \holo\ in a \nbd\ of a \hc\ set $K \subset X$
can be approximated on $K$ by functions holomorphic on $X$ [H\"o2].

Whenever we use a metric on a manifold, we always mean a metric
compatible with the underlying manifold topology. We shall
not mention this again.

The next result is the {\it parametric h-principle} with approximation
on \hc\ sets. A discussion in this direction can be found in
sect.\ 3 in [Gro].

%
%
%
%
\proclaim 1.4 Theorem:  Let $X$ be a Stein manifold and $h\colon Z\to X$ a
locally trivial \holo\ fiber bundle whose fiber admits a spray.
Let $Y$ be a compact Hausdorff space (the parameter space), $Y_0 \subset Y$ a
compact subset and $Y'\subset Y$ an open set containing $Y_0$. Assume that
$f\colon X\times Y\to Z$ is a continuous map such that $f(\cdotp,y) \colon X\to Z$
is a section of $h\colon Z\to X$ for each $y\in Y$, and $f(\cdotp,y)$ is \holo\
on $X$ for each $y\in Y'$. Then there is a homotopy
$H\colon X\times Y\times [0,1] \to Z$ such that
$H_0=f$, $H_1(\cdotp,y)\colon X\to Z$ is \holo\ on $X$ for each $y\in Y$,
and the homotopy is fixed on $Y_0$ (i.e., $H_t(x,y)$ is independent of $t$ for
$y\in Y_0$). Moreover, if $K$ is a compact \hc\ subset in $X$ and we assume that
there is a \nbd\ $V\subset X$ of $K$ such that each section $f(\cdotp,y)$ ($y\in Y$)
is \holo\ in $V$, then for any metric $d$ on $Z$ and for any
$\e>0$ there is a homotopy $H$ as above which also satisfies
$$
    d\bigl( H_t(x,y), f(x,y) \bigr) <
    \e \quad (x\in K,\ y\in Y,\ 0\le t\le 1).
$$

%
%
\ni\it Example 3. \rm  Theorem 1.4 fails for maps of any Stein manifold
$X$ into any Kobayashi-hyperbolic complex manifold $F$. To see this,
take $K$ to be a small piece of an
embedded analytic disc in $X$ and let $f_0 \colon K\to F$ be
a \holo\ map which is close to being extremal
with respect to the Kobayashi metric on $F$. Such $f_0$ cannot be
approximated on $K$ by \holo\ maps $f\colon X\to F$ since this
would create much larger analytic discs in $F$. Of course a
hyperbolic manifold admits no spray since it even admits no
non-degenerate holomorphic images of $\C$.

%
%
\proclaim 1.5 Corollary: If $h\colon Z\to X$ is as in theorem 1.3,
the inclusion map between the spaces of holomorphic and continuous
sections, ${\rm Holo}(X,Z) \hra {\rm Cont}(X,Z)$,
is a {\it weak homotopy equivalence}, i.e., it induces an
isomorphism of the homotopy groups of the two spaces
(endowed with the usual compact--open topology.)

\demo Remarks: 1.  The basic h-principle (def.\ 1.1) is equivalent to saying
that {\it each path connected component of the space of continuous sections
${\rm Cont}(X,Z)$  contains precisely one path connected component of the
space of holomorphic sections ${\rm Holo}(X,Z)$.}

\ni 2.\  If the set $Y_0$ in theorem 1.4 is a deformation retraction
of some open \nbd\ $Y'\subset Y$, it suffices to assume only that
the sections $f(\cdotp,y)$ for $y\in Y_0$ are \holo\ on $X$
since we can use the deformation retraction to suitably reparametrize
the family. In all our applications $Y$ will be a polyhedron and
$Y_0$ a subpolyhedron.

\demo Proof of corollary 1.5:
If we take $Y$ to be the $n$-sphere $S^n$ and $Y'=\emptyset$,
theorem 1.4 implies that each (continuous) map $S^n\to {\rm Cont}(X,Z)$
can be homotopically deformed to a map $S^n\to {\rm Holo}(X,Z)$.
Similarly, if we take $Y$ to be the closed real $(n+1)$-ball
$B^{n+1}$ and $Y_0=\di B^{n+1}=S^n$, we conclude that each map
$S^n\to {\rm Holo}(X,Z)$ which extends to a map
$B^{n+1}\to {\rm Cont}(X,Z)$ also extends to a map
$B^{n+1}\to {\rm Holo}(X,Z)$. This is precisely the content
of corollary 1.5.
\endpr

%
%
%
%

We wish to put these results in historic perspective
as we see it. The term {\it Oka's principle} can vaguely
be stated by saying that,
{\it on a Stein manifold or a reduced Stein space, any analytic
problem with no topological obstruction has an analytic solution.}
To be more precise we quote from p.\ 145 in [GRe]:
{\it Analytic problems which can be cohomologically formulated
have only topological obstructions}. Early examples
include Oka's classification of \holo\ line bundles
and the solvability of the second Cousin problem
[H\"o2], p.\ 144. There are also examples of this principle
on non-Stein manifolds, such as the Riemann-Roch theorem.

A major extension of Oka's principle was obtained by H.\ Grauert
in a series of deep and influential papers [Gr2]--[Gr4] in
1957-58; see also Cartan [Car] for a nice exposition of Grauert's
work. The main result of Grauert [Gr3] was the proof of theorems
1.3 and 1.4 above for holomorphic fiber bundles $E\to X$ over Stein
spaces $X$ (possibly with singularities), where the fiber $L$ is a
complex Lie group and the transition maps $\psi_{\a,\b}(x,\cdotp)$
are Lie group automorphisms of $L$. Each fiber $E_x$
therefore carries a natural Lie group structure isomorphic to $L$,
and there is a well defined identity section of $E$.
The (flows of) left invariant \hvf s on $L$ induce a spray
on $L$, and hence theorem 1.3 includes Grauert's theorem.
In fact there exists even a global spray on $E$, induced by flows
of left invariant fields on $E$ tangent to the fibers $E_x$.

Cartan [Car] observed that Grauert's proof carries over
immediately to fiber bundles $\tilde E\to X$ over a Stein base,
where $E\to X$ is as above and where each fiber $\tilde E_x$
is a $E_x$-homogeneous space, i.e., the elements of $E_x \simeq L$ act
on $\tilde E_x$ by right multiplication. Cartan calls such
$\tilde E$ an {\it espace analytique E-principal.}
(Cartan's observation was acknowledged by Grauert in the sentence
preceeding his Satz 2 in [Gr3] on p.267.)
When $E\simeq X\times L$ is a trivial bundle,
$\tilde E$ is a classical {\it principal $L$-bundle}.
The validity of the h-principle for such bundles $\tilde E$
implies that {\it the holomorphic classification of the principal
holomorphic bundles over a Stein space agrees with the topological
classification} [Car, Gr4]. This holds in particular for holomorphic
vector bundles. Further results were obtained by Ramspott [Ram],
Forster and Ramspott [FR1, FR2], Forster [For], Heinzner and
Kutzschebauch [HKu], and others.

Grauert's constructions were similar in spirit to the usual
proof of Cartan's theorems A and B for coherent analytic sheaves.
The local patching of holomorphic sections on small sets was
obtained (like the `patching of syzygies') by an analogue
of the Cartan's splitting lemma for holomorphic matrices.
The globalization followed the usual scheme of passing from one
compact Stein block in $X$ to the next block. (A Stein block
is a compact set in $X$ that can be embedded as a closed complex
subvariety of a cube in a Euclidean space $\C^N$ by a holomorphic
map $X\to\C^N$.) To insure the convergence one must correct the
section at each step by a suitable version of the Runge approximation
theorem which had been established in [Gr2].

This scheme does not seem to carry over to more general fiber
bundles because no simple proof of the Runge theorem is
available. Things didn't progress very much until 1986
when Henkin and Leiterer ([HL2], [HL3]) introduced the Grauert's
{\it bump method\/} to this problem
and gave a conceptually new proof of Grauert's
theorem. This method had been used earlier in the solution of
the Levi problem [HL1]; we shall comment on it more extensively
below. The crucial point in this approach is that
{\it the Runge theorem is only needed for sections on small
compact subsets of $X$.} Possibly the only shortfall is that,
so far at least, the method only applies to {\it nonsingular}
Stein bases.

In 1989 M.\ Gromov [Gro] made a crucial step by replacing
the exponential map in fibers by the much more flexible concept
of a spray. This made it possible to extend the results to
the present form.

When reading sect.\ 2 in [Gro] we were unable to complete the proof
as suggested there due to seemingly nontrivial analytic and
geometric problems. Eventually we completed the proof
in a different way, using the work of Henkin and Leiterer [HL3]
and the tools from [Gro]. We give here a complete exposition of
our version of the proof. We also take this opportunity to provide
a detailed account, with proofs, of all major tools.

In the rest of this section we explain the outline of our proof
and also point out the difficulties that we had with
sect.\ 2 in [Gro]. For the sake of simplicity we only discuss
the simplest case: to deform a continuous section by a homotopy to
a holomorphic section.

The underlying geometric scheme is the Grauert's `bump method'
(sect.\ 2); our reference for this is [HL3]. The manifold $X$ is
exhausted by an increasing family of smooth, compact, \spsc\ domains
$A_k \ss X$ such that $A_{k+1}=A_k\cup B_k$,
where $B_k$ is a small \spsc\ domain
attached to $A_k$ in a certain special way (a {\it special \psc\ bump}
in the terminology of [HL3]; see def.\ 2.6 there or def.\ 2.2 in
the present paper). We may assume that the bundle $Z$ is trivial
over a \nbd\ of $B_k$. To obtain such an exhaustion one begins
with a smooth \spsh\ exhaustion function $\rho\colon X\to \R$
with nice critical points. One can pass from one sublevel set
$\{\rho \le c\}$ to a higher sublevel set $\{\rho\le c'\}$,
where $c<c'$ are regular values of $\rho$, by attaching a small
`bump' $B_k$ at each step. These bumps are contractible
\spsc\ domains such that $C_k=A_k\cap B_k$ is also contractible, except
when passing a critical point of $\rho$ when $C_k$ becomes a torus.
Moreover, $C_k$ is Runge in $B_k$.

Granted such an exhaustion of $X$,
we inductively construct a sequence
of continuous sections $f_k\colon X\to Z$
such that $f_k$ is holomorphic in a \nbd\ of $A_k$, it approximates
$f_{k-1}$ on $A_{k-1}$, and is homotopic to the original section $f_0$.
The desired holomorphic section is then obtained as the locally uniform
limit $f=\lim_{k\to\infty} f_k \colon X\to Z$.

Suppose $f_k$ as above has been constructed. To construct
$f_{k+1}$ we proceed as in [Gro] or [HL3]. Since $Z$ is trivial
over $B_k$, there is a holomorphic section $b_0\colon B_k\to Z$
homotopic to $f_0$. The first step is to construct a
{\it holomorphic} homotopy $b_t\colon C_k\to Z$, $0\le t\le 1$,
connecting $b_0$ and $b_1=f_k|C_k$ over $C_k$. This is the
difficult part and it is here that our proof differs from
[Gro]. We shall comment on this later. From now on the process
follows [Gro]. Since $b_0$ is holomorphic on $B_k$, and $C_k$
is Runge in $B_k$, the homotopy version of Runge theorem (sect.\ 4)
allows us to approximate the homotopy $b_t$ uniformly on $C_k$ by
a holomorphic homotopy $\tilde b_t$ defined on $B_k$. The
section $b=\tilde b_1$ approximates $f_k$ as well as desired on
$C_k$, so we can glue them (see sect.\ 5) into a single section
$f_{k+1}$ that is holomorphic on $A_{k+1}=A_k\cup B_k$ and approximates
$f_k$ on $A_k$. It remains to extend $f_{k+1}$ continuously
to $X$ so that it is homotopic to $f_0$. In the process we must
also insure the convergence of homotopies from $f_0$ to $f_k$
on $X$.

Back to the homotopy $b_t$ on $C_k$. In the {\it non-critical case}
the set $C_k=A_k\cap B_k$ is holomorphically contractible.
Since $Z$ is trivial over $C_k$, we immediately get the required
homotopy by using the contractions of $C_k$ to a point.
This non-critical case allows us to proceed in a finite number
of steps from a sublevel set $\{\rho\le c_0\}$ to a higher
sublevel set $\{\rho \le c_1\}$, provided that $\rho$ has no
critical values on the interval $[c_0,c_1] \subset \R$.
A similar method (with $C_k=\emptyset$) allows us to cross
the critical points of $\rho$ which are local minima.

To pass a critical point $x_0$ of $\rho$ which is not a local
minimum, Gromov suggested (sec.\ 2.7 in [Gro]) to attach to
the set $A_k = \{\rho\le c\}$, for some $c<c_0=\rho(x_0)$
and $c$ close to $c_0$, a \ra, \tr\ disc $B_k$ containing $x_0$.
We must then find small \spsc\ \nbd s $V \subset X$ of $A_k\cup B_k$
such that we can reach a suitable higher sublevel set $\{\rho\le c'\}$
for $c'> c_0$ by attaching convex bumps to $V$ (i.e., as a
non-critical pseudoconvex extension of $V$).
If we understand correctly, the idea in [Gro]
was to take a homotopy $\{b'_t\}$ from $b'_0=b_0$ to
$b'_1=f_k$ over the disc $B_k$ (such exist since $B_k$
is contractible), then approximate it by a \ra\ homotopy
$b_t$ on $B_k$, and finally complexify $b_t$
to get a holomorphic homotopy whose final section $b=b_1$
approximates $f_k$ on (and hence near) the \tr\ collar
$C_k=A_k\cap B_k$. It was then suggested to proceed
with gluing as before.

We could not complete this for two reasons. The first
problem is to find suitable \nbd s $V$ of $A_k\cup B_k$
with the required properties. We do not know of any such result
in the literature; for partial results see [Eli] and [Ros].
Another independent problem is that the rate of approximation of
$f_k$ by $b$, required to glue them into a single section,
depends on $V$ since we must solve a certain $\dibar$-equation
with sup norm estimates in the process. Even though this can be
done on any \spsc\ domain, the constant in such estimate depends
(unlike for the $L^2$ estimate!) on the geometry of the set.
If the sets degenerate, the constants will blow up in general.
If on the other hand we work with the $L^2$ estimates, the
loss in the Cauchy estimates again depends on the shape of $V$.
However, we do not see how to control the rate of approximation
of $f_k$ by $b$ in any fixed \nbd\ of the collar $A_k\cap B_k$
when $b$ is obtained as above. See also remark 1 following the
proof of lemma 2.4 below.

We instead proceed as follows. Suppose that $x_0 \in X$ is a critical
point of $\rho$. Set $c_0=\rho(x_0)$. It was proved in [HL3]
that, if $c<c_0$ and $c$ is sufficiently close to $c_0$,
we can attach to the set $A_k =\{x\in X\colon \rho(x)\le c\}$
a  \psc\ bump $B_k$ as follows. The set $C_k=A_k\cap B_k$ is a
sublevel set $\{\tau\le 1\}$ of some \spsh\ function $\tau\ge 0$,
defined in a \nbd\ of $B_k$, such that (in some holomorphic coordinate
system in a \nbd\ of $B_k$) the set $S=\tau^{-1}(0)$ is a
sphere contained in an affine \tr\ subspace, and $\tau$ has no
critical points on $C_k\bs S$. In particular, $C_k$ is obtained
from a small tubular \nbd\ of $S$ by attaching convex bumps.
Furthermore, a suitable sublevel set $\{\rho\le c'\}$ for some
$c'>c_0$ is a non-critical \spsc\ extension of $A_{k+1}=A_k\cup B_k$.

To construct the homotopy $b_t$ over $C_k$ we initially
deform a given continuous homotopy $b'_t$ from $b'_0=b_0$
to $b'_1=f_k$ (which exists on $B_k$) in a small tube around
the sphere $S$ to make it \holo\ there. This is possible
since $S$ is \tr. Subsequently we extend it to $C_k$ in a finite
number of steps. In each step we extend the homotopy
(by approximation) across a convex bump, using parametric
versions of the h-principle and of the gluing lemma.
Since $C_k$ is a non-critical extension of a tube around $S$,
this allows us to complete the proof of theorems 1.3 and 1.4.
The crucial new steps in our proof are theorem 4.5 and corollary 5.6.

The paper is organized as  follows. In section 2 we recall the
{\it bump method}, following [HL3], and we reduce the proof of
theorem 1.4 to theorem 2.6 concerning the extension
of holomorphic sections across \pb s.
Theorem 2.6 is proved in sect.\ 6 after we develop the necessary
tools in sect.\ 3--5. In sect.\ 3 we recall from [Gro] the
relevant properties of sprays and iterated spray bundles,
and we prove that any homotopy of \holo\ sections of $Z\to X$ can be
lifted to a homotopy of sections of an iterated spray bundle over $Z$.
In sect.\ 4 we prove Runge-type approximation theorems for sections
of holomorphic submersions which admit a spray. In sect.\ 5 we
prove results on gluing holomorphic sections over Cartan pairs.
In sect.\ 6 we complete the proof of theorem 2.6.

We wish to thank G.\ Henkin and J.\ Leiterer for their interest in
our work and for several useful discussions. We also thank J.\ Globevnik
who got us interested in this topic, and the participants in the
Seminar for complex analysis at the University of Ljubljana who
suffered through our lectures. The first author acknowledges partial
support by the National Science Foundation, by the Vilas Foundation at
the University of Wisconsin, and by the Ministry of Science and Technology
of the Republic of Slovenia.  The second author was supported by the
Ministry of Education of the Republic of Slovenia.

%
%
%
%
\beginsection 2. Pseudoconvex bumps and Cartan pairs.

The main reference for this section is [HL3]. Let $X$ be
a complex manifold. We say that a compact set $C\subset X$ is Runge in
another compact set $B\supset C$ if $C$ has a basis of \nbd s $C_j$ which are
Runge in some open \nbd\ $\tilde B$ of $B$, i.e., each function holomorphic in
$C_j$ can be approximated uniformly on compacts in $C_j$ by functions holomorphic
in $\tilde B$. This is the case for instance if $C$ is \hc\ in some Stein \nbd\
of $B$.

\medskip\ni \bf  2.1 Definition. \sl Let $X$ be a complex manifold.

\item{(i)} A {\bf compact \spsc\ domain} in $X$ is a compact set of the form
$A=\{\rho\le 0\} \subset X$, where $\rho\colon X\to \R$ is a $\cC^2$ function
which is \spsh\ in a \nbd\ of $\{\rho=0\}$ and has no critical points
on $\{\rho=0\}$.

\item{(ii)}
Let $A\subset A'$ be a pair of compact \spsc\ domains in $X$.
We say that $A'$ is a {\bf \spsc\ extension} of $A$ in $X$
if there is a $\cC^2$ function $\rho\colon X\to\R$
which is \spsh\ on $\bar {A'\bs A}$ such that, for some
real numbers $\beta > \alpha$ which are regular values
of $\rho$,
$$  A= \{x\in X\colon \rho(x)\le \alpha \},\quad
    A'=\{x\in X \colon \rho(x)\le \beta \}.
$$
If $\rho$ can be chosen such that it has no critical
points on $\bar{A'\bs A}$, $A'$ is called
a {\bf non-critical \spsc\ extension} of $A$.

%
%
\medskip\ni \bf 2.2 Definition. \sl  Let $X$ be a complex manifold.
A pair $(A,B)$ of compact subsets of $X$ is a {\bf \psc\ bump}
(or $B$ is a \psc\ bump on $A$) if the following hold:
\item{(i)}    The sets $A$, $B$, $A\cup B$ and $C=A\cap B$ are compact \spsc\
domains ($C$ may be empty);
\item{(ii)}   $\bar{A\bs B} \cap \bar{B\bs A} =\emptyset$;
\item{(iii)}  there are \holo\ coordinates in a \nbd\ of $B$ in which
$B$ is star-shaped;
\item{(iv)}   $C$ is Runge in $B$;
\item{(v)}  if $C \ne\emptyset$ there is a $\cC^2$ \spsh\ function $\tau\ge 0$
defined in a \nbd\ $U\subset X$ of $B$ such that $C=\{z\in U\colon \tau(z) \le 1\}$,
the set $S=\tau^{-1}(0)$ is a compact \tr\ submanifold contained in an affine
\tr\ subspace (with respect to some holomorphic coordinates in a \nbd\ of $B$),
and $\tau$ has no critical points in $C\bs S=\{z\in U\colon 0<\tau(z)\le 1\}$.

\ni We say that $B$ is a {\bf convex bump} on $A$ (or the pair $(A,B)$ is a convex bump)
if, in addition to the above, there are holomorphic coordinates in a \nbd\ of $B$
in which both $B$ and $C$ are strongly convex domains. (The set $C$ may be empty.)
\medskip\rm

For certain purposes we can relax the above conditions and consider
{\it Cartan pairs}. Gromov's definition of a Cartan pair is rather
imprecise (sec.\ 1.5.A in [Gro]); we shall adopt the following
definition.

%
%
%
%
\proclaim 2.3 Definition:
A {\bf Cartan pair} in $X$ is a pair of compact sets $A,B\subset X$ such that
\item{(i)}    $A$, $B$, and $A\cup B$ have bases of Stein \nbd s, and
\item{(ii)}   $\bar{A\bs B} \cap \bar{B\bs A} =\emptyset$.

Any \psc\ bump is clearly a Cartan pair. We denote by $H^\infty(\Omega)$ the
algebra of bounded \holo\ functions on $\Omega$. The following lemma will
be used in sec.\ 5 for gluing \holo\ sections over Cartan pairs;
this is similar to what Gromov takes as the definition of a Cartan pair.

\proclaim 2.4 Lemma: If $(A,B)$ is a Cartan pair in a Stein manifold
such that $C=A\cap B\ne\emptyset$, there are
bases of Stein open \nbd s $A_j\supset A$, $B_j\supset B$,
$C_j =A_j\cap B_j \supset A\cap B$, and bounded linear operators
${\cal A}_j\colon H^\infty(C_j) \to H^\infty(A_j)$,
${\cal B}_j\colon H^\infty(C_j) \to H^\infty(B_j)$,
satisfying
$$ c= {\cal A}_j(c) - {\cal B}_j(c),\qquad c\in H^\infty(C_j),\ j=1,2,3,\ldots
                                                                   \eqno(2.1)
$$

\demo Proof: Let $U\supset A$ and $V\supset B$ be open \nbd s of $A$ resp.\ $B$.
Choose Stein open sets  $\tilde A$ and $\tilde B$ in $X$ so that
$A\subset \tilde A\subset U$ and $B\subset \tilde B\subset V$. Set
$\tilde C=\tilde A\cap \tilde B$. By the separation condition (ii) in def.\ 2.3
there is a smooth function $\chi\colon X\to [0,1]$ such that $\chi=0$ in a \nbd\
of $\bar{A\bs B}$ and $\chi=1$ in a \nbd\ of $\bar{B\bs A}$. Hence there
are open sets $A_0, B_0\subset X$, with $A\subset A_0\subset \tilde A$ and
$B\subset B_0\subset\tilde B$,
such that $\chi=0$ on $A_0\bs \tilde B$ and $\chi=1$ on
$B_0\bs\tilde A$. Choose a smooth \spsc\ set $\Omega$, with
$A\cup B\subset \Omega\subset A_0\cup B_0$, and set
$$ A'=\tilde A\cap \Omega,\quad B'=\tilde B\cap\Omega,\quad
   C'=A'\cap B'= \tilde C\cap \Omega.
$$
Then $A'$, $B'$ and $C'$ are Stein domains containing $A$, $B$, $C$
respectively, with $A'\subset \tilde A\subset U$,
$B'\subset \tilde B\subset V$, and
$A'\cup B'=\Omega$. Moreover by the choice of $\Omega$ we have
$$ A'\bs C'= A'\bs B'=
   \Omega\bs\tilde B \subset A_0\bs \tilde B \subset \{\chi=0\}
$$
and similarly $B'\bs C'\subset \{\chi=1\}$.
If $c$ is any bounded \holo\ function
in $C'$ then the above implies that $\chi c$ extends to
a bounded smooth function in $A'$ which vanishes in
$A'\bs C'$, $(1-\chi)c$ extends to a bounded smooth function
in $B'$ which vanishes in $B'\bs C'$, and $\dibar(\chi c)=c\dibar\chi$
extends to a bounded smooth $\dibar$-closed $(0,1)$-form in $\Omega$.
Let $T_\Omega$ be a linear solution operator for the $\dibar$-equation
in $\Omega$ (i.e., $\dibar(T_\Omega \alpha)=\alpha$ for each
$\dibar$-closed smooth $(0,1)$-form $\alpha$ on $\Omega$) which is
bounded in the sup norm (and also in any $\cC^k$ norm);
see [HL1], p.\ 82 or [H\"o1]. The linear operators
$$  {\cal A}\,c = \chi c    - T_\Omega(c\,\dibar \chi),\qquad
    {\cal B}\,c = (\chi-1)c - T_\Omega(c\,\dibar \chi)
$$
then satisfy the required properties with respect to the \nbd s $A'$, $B'$, $C'$.
In fact, since $T_\Omega$ is bounded in all $\cC^k$ norms, so are the operators
$\cal A$ and $\cal B$.
\endpr

\demo Remarks: 1.\ Gromov requires the existence of decompositions
(2.1) with the sup norm estimates on ${\cal A}_j(c)$ and ${\cal B}_j(c)$
only depending on $A$ and $B$ and not on the neighborhoods
(see sect.\ 1.5.A in [Gro]). The proof of lemma 2.4 shows that
this is possible if we can solve the $\dibar$-equation with
the sup norm estimates on a system of \nbd s of $A\cup B$
so that the constant in these estimates only depends on $A$ and
$B$. This can be achieved only rarely, for instance when $A\cup B$
is itself a \spsc\ domain. This problem causes difficulties
in  sect.\ 2.7 of [Gro] when crossing a critical point.
The remark 1.5.A' in [Gro] is incorrect, with immediate counterexamples.
(The union of two \hc\ sets in $X$ need not have a basis of Stein
\nbd s, and there may be no decomposition (2.1).)

\ni 2.\ Instead of the Banach space $H^\infty$ we could
as well use in lemma 2.4 any of the spaces $A^k$
(holomorphic functions which are smooth of class $\cC^k$
up to the boundary) as in [HL3]. The solution operator
$T_\Omega$ on any \spsc\ domain is bounded in any of
these norms.

The inductive construction of sections of $Z\to X$ is based on the
following geometric result from [HL3].

\proclaim 2.5 Theorem:
Let $X$ be a complex manifold and let $A'\subset X$ be a \spsc\
extension of a subset $A\subset A'$ (def.\ 2.1 (ii)).
Then there exist \psc\ bumps $(A_j,B_j)$ in $X$, $0\le j\le k$,
such that
$$ A=A_0,\quad A_{j+1}=A_j\cup B_j\ \ {\rm for\ \ } 0\le j\le k-1,\quad
   {\rm and \ \ } A'=A_k.
$$
Moreover for any open covering $\{U_i\}$ of $X$ we can choose the bumps such that
each $B_j$ is contained in some $U_i$.  If $A'$ is a non-critical \spsc\ extension
of $A$ then we may choose each $B_j$ to be a convex bump on $A_j$.

The following is the main technical result of this paper. It allows us to
extend a holomorphic section of $Z\to X$ from a \nbd\ of $A$ to a \nbd\ of
$A\cup B$ for each \psc\ bump $(A,B)$ in $X$, with approximation on $A$.

\proclaim 2.6 Theorem:
Let $h\colon Z\to X$ be a holomorphic submersion onto a Stein manifold $X$ and let
$(A,B)$ be a \psc\ bump in $X$. Assume that there is an open \nbd\ $\tilde B \subset X$
of $B$ such that $Z|\tilde B=h^{-1}(\tilde B)$ is isomorphic to a trivial bundle
$\tilde B\times F$, where $F$ admits a spray. Let $Y$ be a compact Hausdorff
space (the parameter space), $Y_0\subset Y$ a compact subset,
and $Y'\subset Y_0$ an open set containing $Y_0$. Let $U\subset X$ be a
\nbd\ of $A$. Suppose that $a\colon X\times Y\to Z$ is a continuous map such
that for each $y\in Y$, $a(\cdotp,y) \colon X\to Z$ is a section of $h\colon Z\to X$
which is \holo\ in $U$, and the sections $a_y$ for $y\in Y'$ are \holo\ on $X$.
Then for any metric $d$ on $Z$ and for any $\e>0$ there exists a
homotopy $a_t\colon X \times Y\to Z$ ($0\le t\le 1$) satisfying
\item{(i)}   $a_t(\cdotp,y)\colon X\to Z$ is a section of $h\colon Z\to X$
for each $y\in Y$ and $0\le t\le 1$,
\item{(ii)}  $a_0=a$,
\item{(iii)} each section $a_1(\cdotp,y)$ is \holo\ in a \nbd\ of $A\cup B$
(independent of $y$),
\item{(iv)}  the homotopy is fixed for $y$ in a \nbd\ of $Y_0$, i.e.,
$a_t(\cdotp,y)=a(\cdotp,y)$ for $y$ near $Y_0$ and $0\le t\le 1$, and
\item{(v)}   $d\bigl( a_t(x,y), a(x,y) \bigr) < \e$ for $x\in A$, $y\in Y$ and
$0\le t\le 1$.
\medskip\rm

Theorem 2.6 is proved in sect.\ 6 below. In the rest of this section we
assume that theorem 2.6 holds and prove theorem 1.4. For simplicity
we write the proof in the case without the parameter $y$;
the proof of the general case is the same.

\demo Proof of theorem 1.4:
We may assume that $d$ is complete metric on $Z$. Fix an open
covering ${\cal U}=\{U_i\}$ of $X$ such that $Z|U_i$ is a trivial
bundle with fiber $F$ for each $i$. Since
$X$ is Stein and $K$ is \hc\ in $X$, there is a smooth \spsh\ exhaustion function
$\rho\colon X\to \R$ such that $\rho<0$ on $K$, $0$ is a regular value of $\rho$, and
the given section $f_0$ is holomorphic in a \nbd\ of
$A_0=\{x\in X\colon \rho(x)\le 0\}$. By theorem 2.5 there is a sequence of
compact \spsc\ domains
$A_0\subset A_1\subset A_2\subset\ldots \subset \bigcup_{k=0}^\infty A_k  = X$
such that $A_{k+1}=A_k\cup B_k$ for each $k$, where $(A_k,B_k)$ is a \psc\ bump
and $B_k\subset U_i$ for some $i$ (so that $Z$ is trivial over a \nbd\ of $B_k$).

Applying theorem 2.6 to each \pb\ $(A_k,B_k)$ we can inductively construct a
sequence of continuous sections $f_k\colon X\to Z$ (with $f_0$ being the given
initial section) and homotopies of continuous sections
$H^k_t\colon X\to Z$ ($0\le t\le 1$) satisfying the following properties
for each $k=0,1,2,\ldots$ and $0\le t\le 1$:

\medskip
\item{(i)}    $f_k$ and $H^k_t$ are holomorphic in a \nbd\ of $A_k$;
\item{(ii)}   $H^k_0=f_k$, $H^k_1=f_{k+1}$;
\item{(iii)}  $d\bigl( H^k_t(x), f_k(x) \bigr) < \e/2^{k+1}$ for all $x\in A_k$
and $0\le t\le 1$.
\medskip

We begin by applying theorem 2.6 to the data $a=f_0$, $A=A_0$, $B=B_0$, and
$\e$ replaced by $\e/2$ to get a homotopy $H^0_t$ such that the section
$H^0_1=f_1$ is holomorphic in a \nbd\ of $A_0\cup B_0=A_1$ and the properties
(i)--(iii) hold with $k=0$.

Suppose inductively that $k>0$ and we have already constructed sections
$f_l$ and $H^l_t$ for $l<k$ and $0\le t\le 1$. Then we take $f_k=H^{k-1}_1$
(this section is holomorphic in a \nbd\ of $A_k$) and apply theorem 2.6
to the data $a=f_k$, $A=A_k$, $B=B_k$, and $\e$ replaced by
$\e/2^{k+1}$, to get $H^k_t$ and $f_{k+1}=H^k_1$ satisfying (i)--(iii).
This completes the induction step.

Property (iii) with $t=1$ implies for all $k=0,1,2,\ldots$
$$ d\bigl( f_{k+1}(x), f_k(x) \bigr) < \e/2^{k+1} \quad (x\in A_k). $$
Since the metric $d$ is complete, it follows that the limit \
$f=\lim_{k\to\infty} f_k \colon X\to Z$\ exists uniformly on compacts in $X$,
it is holomorphic on $X$ since $f_k$ is holomorphic on $A_l$ for each $k\ge l$,
and it satisfies the estimate $d\bigl( f(x), f_0(x) \bigr) <\e$ for $x\in A_0$.

To construct a homotopy $H_t$ between $H_0=f_0$ and $H_1=f$ we divide
the $t$-interval $[0,1)$ into subintervals
$I_k=[1-2^{-k},1-2^{-k-1}]$ ($k=0,1,2,\ldots$) and take $H_t$ for
$t\in I_k$ to be the homotopy $H^k_t$, suitably rescaled. To be precise, let
$\l_k \colon I_k\to [0,1]$ be the linear bijection
$\l_k(t)=2^{k+1}\bigl( t- 1+2^{-k} \bigr)$.
For each $k=0,1,2,\ldots$ we set
$$ H_t(x)= H^k_{\l_k(t)} (x) \quad   (t\in I_k,\ x\in X). $$
Clearly this defines a homotopy $H_t\colon X\to Z$ for $0\le t<1$.
The property (iii) above implies that $\lim_{t\to 1} H_t =f$ uniformly on
compacts in $X$, so by setting $H_1=f$ we obtain the required homotopy
from $f_0$ to $f$. This completes the proof of  theorem 1.4, provided
that theorem 2.6 holds.
\endpr

%
%
%
%
\beginsection 3. Holomorphic submersions with sprays.

The reference for this section are sections 1.1-1.3 in [Gro].
We fist recall the notion of a (fiber dominating) spray associated
to a holomorphic submersion $h\colon Z\to X$. We also recall
the notion of an {\it iterated spray bundle} and we prove results
on lifting homotopies to iterated spray bundles. This is used
in sect.\ 4 in the proof of h-Runge theorems. While it would be
possible to avoid iterated spray bundles, their use simplifies
proof of theorem 4.2.

Let  $h\colon Z \to X$ be a \holo\ submersion between complex manifolds
(not necessarily Stein). For $x\in X$ we denote by $Z_x=h^{-1}(x) \subset Z$
the fiber over $x$.  At each point $z\in Z$ the
tangent space $T_z Z$ contains a well defined {\bf vertical tangent space}
$$  VT_z(Z) = \{ e\in T_z Z\colon Dh(z)e =0\} = T_z Z_{h(z)}.           \eqno(3.1) $$
We denote by $VT(Z)$ the corresponding {\bf vertical tangent bundle} to $Z$
which is a \holo\ subbundle of the tangent bundle $TZ$. Since $Z$ is not
assumed to be Stein, there is in general no splitting of $TZ$ into a direct
sum $VT(Z)\oplus E'$ for some \hvb\ $E'\mapsto Z$. However such a splitting exists
over any open Stein subset (or a Stein submanifold) $V\subset Z$ (see [GRo], p.256).
Also if $f\colon X\to Z$ is a holomorphic section then along the graph $f(X)\subset Z$
the tangent bundle has a canonical splitting
$$
     TZ|f(X) = VT(Z)|f(X) \oplus Tf(X).
$$
If $p\colon E\to Z$ is a \hvb\ over $Z$, we denote by $E_z=p^{-1}(z) \subset E$
its fiber over $z\in Z$ and by $0_z \in E_z$ the zero element of $E_z$.

%
%
\medskip\ni \bf 3.1 Definition.  \sl
A {\bf spray} on $Z$ associated to the submersion $h\colon Z\to X$
(or a {\bf fiber-spray}) is a tripple $(E,p,s)$, where $p\colon E\to Z$ is a \hvb\ and
$s\colon E\to Z$ is a \holo\ map such that for each $z\in Z$ we have
\item{(i)}   $s(E_z) \subset Z_{h(z)}$ (equivalently, $h\circ p=h\circ s$),
\item{(ii)}  $s(0_z)=z$, and
\item{(iii)} the derivative $Ds(0_z) \colon T_{0_z} E \to VT_z(Z)$
maps the subspace $E_z \subset T_{0_z}E$ surjectively onto $VT_z(Z)$.
\medskip\rm

%
%
\ni We denote the restriction in (iii) by
$$
    VD s(z) = Ds(0_z)|E_z \colon E_z\to VT_z(Z)                         \eqno(3.2)
$$
and call it the {\bf vertical derivative} of $s$ at the point $0_z\in E$.
Gromov [Gro] calls such a map $s$ a {\it fiberwise dominating spray}, the word
dominating referring to the property (iii). We shall call it simply a spray
when there is no danger of confusion with def.\ 1.2, or a fiber-spray
if we wish to emphasize the difference between the two notions.

\demo Example: Each spray $s \colon E\to F$ in the sense of definition 1.2
induces a fiber-spray $(\tilde E,\tilde p,\tilde s)$ associated to the
trivial fibration
$h \colon Z=X\times F\to X$ by taking
$$ \tilde E= X\times E, \quad \tilde p(x,e)=(x,p(e))\in X\times F=Z,\quad
   \tilde s(x,e)=(x,s(e)) \in Z.
$$
Hence if $Z\mapsto X$ is a locally trivial bundle whose fiber admits
a spray, then $X$ can be covered by open sets $U_i$ such that each restriction
$Z|U_i$ admits a (fiber-) spray (but in general there is no global spray over $Z$).
\endpr

The main use of sprays is to lift homotopies of sections of $h\colon Z\to X$
to homotopies of sections of a certain vector bundles, thereby linearizing
the approximation and gluing problems for such sections.
The first result in this direction is

\proclaim 3.2 Lemma: \rm (Gromov [Gro], sec.\ 1.2.) \sl
Let $X$ be a Stein manifold and $h\colon Z\to X$ a \holo\ submersion which
admits a spray $(E,p,s)$. Then for each \holo\ section $f\colon X\to Z$
there exists a \holo\ vector subbundle $E'$ of the restricted bundle
$E|f(X)$ such that $s\colon E'\to Z$ maps a \nbd\ of the zero section
in $E'$ biholomorphically onto a \nbd\ of $f(X)$ in $Z$. In particular,
if $f_t\colon X\to Z$ ($0\le t\le 1$) is a homotopy of holomorphic
sections, then for each $t_0\in [0,1]$ and each open relatively compact
subset $V\ss X$ there is \nbd\ $I_0\subset [0,1]$ of $t_0$ and a homotopy
of holomorphic sections $\xi_t$ ($t\in I_0$) of $E'\subset E|f_{t_0}(X)$
over the set $f_{t_0}(V)$ such that $\xi_{t_0}$ is the zero section and
$s\circ \xi_t(z)=f_t(h(z))$ for $t\in I_0$ and $z\in f_{t_0}(V)$.

\demo Proof: By definition of the spray the map $s\colon E|f(X)\to Z$
is the identity on the zero section (which we identify with $f(X)\subset Z$)
and it is a submersion near the zero section.
Denote by $E_0=\ker VDs \subset E$ the kernel of the vertical
derivative (3.2)  and let $\tilde E=E/E_0$ be the
quotient bundle with the quotient projection $\pi\colon E\to \tilde E$.
Since $X$ is Stein, this projection splits over $f(X)$, i.e.,
there is a holomorphic vector bundle homomorphism
$G\colon \tilde E|f(X)\to E|f(X)$ such that $\pi\circ G$ is the
identity on $\tilde E|f(X)$. If we denote by $E'$
the image of $G$, we have a direct sum decomposition
$$           E|f(X) = E_0|f(X) \oplus E'.              \eqno(3.3) $$
The restriction $s|E' \colon E' \to Z$ maps the zero section of
$E'$ onto $f(X)$ and its derivative is an isomorphism at each point of the
zero section. Hence $s|E'$ is biholomorphic near the zero section.
The second statement follows immediately from this.
\endpr

Lemma 3.2 allows us to lift short pieces of a homotopy of sections of $Z$
to a homotopy of sections of a vector bundle. In order to lift the entire
homotopy we recall from [Gro] the concept of composed and iterated sprays.

\proclaim 3.3 Definition: \rm (Gromov [Gro], sec.\ 1.3.) \sl
(a) Let $(E_1,p_1,s_1)$ and $(E_2,p_2,s_2)$  be sprays on $Z$ associated to a
submersion $h\colon Z\to X$. The {\bf composed spray} $(E^*,p^*,s^*)$ over $Z$
is defined by
$$ \eqalign{ & E^* = \{(e_1,e_2)\in E_1\times E_2 \colon\ \ s_1(e_1)=p_2(e_2) \}, \cr
             & \ \ p^*(e_1,e_2) = p_1(e_1), \quad s^*(e_1,e_2) = s_2(e_2). \cr}
$$
(b) Let $(E,p,s)$ be a spray on $Z$ associated to $h\colon Z\to X$.
For each integer $k=1,2,3,\ldots$ the $k$-th {\bf iterated spray}
$(E^{(k)},p^{(k)},s^{(k)})$ is defined by
$$ \eqalignno{ & E^{(k)} = \{e=(e_1,e_2,\ldots,e_k) \colon e_j\in E\
            \ {\rm for}\ j=1,2,\ldots,k, \cr
            &\qquad \qquad\qquad
              s(e_j)=p(e_{j+1})\ {\rm for}\ j=1,2,\ldots,k-1\},  \cr
            & p^{(k)}(e) = p(e_1),\quad s^{(k)}(e) = s(e_k). &(3.4) \cr}
$$
\rm

Note that the composed spray is {\it not a spray} over $Z$ in the sense of def.\ 3.1
because $E^*$ does not have a natural structure of a \hvb\ over $Z$ with respect to
the projection $p^*\colon E^*\to Z$ (the other requirements are satisfied). In fact
$E^*$ is the pullback of the vector bundle $p_2\colon E_2\to Z$ by the spray map
$s_1\colon E_1\to Z$, so it is a \hvb\ over $E_1$ with the projection
$(e_1,e_2)\mapsto e_1$. Similarly we can define the iterated sprays inductively
as the composition of $k$ copies of $(E,p,s)$. We begin by taking
$(E^{(1)},p^{(1)},s^{(1)})=(E,p,s)$. Suppose that $(E^{(k-1)},p^{(k-1)},s^{(k-1)})$
has already been defined. Let $q^{(k)} \colon E^{(k)} \to E^{(k-1)}$ be the pullback
of the bundle $p\colon E\to Z$ by the spray map $s^{(k-1)} \colon E^{(k-1)}\to Z$.
Set $p^{(k)}=p^{(k-1)}\circ q^{(k)}  \colon E^{(k)} \to Z$  and let
$s^{k}\colon E^{(k)}\to Z$ be the map induced by $s\colon E\to Z$ under the pullback.
This gives the next iterate $(E^{(k)},p^{(k)},s^{(k)})$.

The following lemma implies that the restriction of composed and iterated
sprays to Stein subsets of $Z$ admit a \hvb\ structure.

\proclaim 3.4 Lemma: Let $Y$ be a Stein manifold and let $p_1\colon E_1\to Y$
resp.\ $p\colon E\to E_1$ be \hvb s over $Y$ resp.\ over $E_1$. Then $E$
has the structure of a \hvb\ over $Y$ with respect to the projection
$p_1\circ p\colon E\to Y$. In fact, this bundle is isomorphic to the
Whitney sum $E_1\oplus (E|Y)$, where $E|Y$ denotes the restriction
of $E$ to the zero section $Y\subset E_1$ of $E_1$.

%
%
\proclaim 3.5 Corollary: {\rm (Gromov [Gro], sec.\ 1.3.A')}
The restriction of any composed or iterated spray bundle on $Z$
to any Stein subset $U\subset Z$ admits a structure of a \hvb\
over $U$.

\demo Proof of lemma 3.4:
Since $E_1\to Y$ is a \hvb\ over a Stein manifold $Y$, the total space
$E_1$ is itself a Stein manifold. Denote the points of $E_1$ by
$(y,e)$, where $y\in Y$ and $p_1(y,e)=y$. Let
$h_t\colon E_1\to E_1$ for $t\in \C$ be the homotopy
$h_t(y,e)=(y,te)$. Consider the family of pull-backs
$h_t^*(E) \to E_1$ of the vector bundle $p\colon E \to E_1$.
It follows from Grauert's theorem [Gr4] (Satz I) that all bundles
in the family are holomorphically isomorphic; see [Lei]
for an elementary proof. In particular the bundles $h_1^*(E)$
and $h_0^*(E)$ are isomorphic. The map $h_1$ is the identity on
$E_1$ whence $h_1^* E=E$. The map $h_0=p_1$ is the projection of
$E_1$ onto the zero section $Y\subset E_1$ and hence
$h_0^*(E)=p_1^*(E)=p_1^*(E|Y)$. The latter bundle
is clearly isomorphic to the Whitney sum $E_1 \oplus (E|Y)$
which is a \hvb\ over $Y$. This proves lemma 3.4.
\endpr

The next result indicates the main application of iterated sprays.
%
%
%
%
\proclaim 3.6 Proposition: Let $f_t\colon X\to Z$ $(t\in [0,1])$
be a homotopy of \holo\ sections of a \holo\ submersion $h\colon Z\to X$.
Assume that the base $X$ is Stein and that $h\colon Z\to X$ admits
a spray $(E,p,s)$. Then for each open relatively compact subset
$V\ss X$ there are an integer $k>0$ and a homotopy of \holo\
sections $\xi_t$ ($0\le t \le 1$) of the iterated spray bundle
$E^{(k)}$ (3.4) over the set $f_0(V) \subset Z$ such that
$$ \xi_0(z)= z,\quad s^{(k)}(\xi_t(z))=f_t(h(z))\qquad
   (z\in f_0(V),\ 0\le t\le 1).
$$

\demo Proof:
For each fixed $t\in [0,1]$ we can apply lemma 3.2 to lift the sections $f_\tau$
for $\tau$ near $t$ by the spray map $s$ to a homotopy of \holo\ sections of
$E|f_t(V)$. Hence by compactness of $[0,1]$ there are numbers
$0=t_0<t_1<t_2<\cdots<t_k=1$ such that for each $j=0,1,\ldots,k-1$ there exists a
homotopy of holomorphic sections $\xi^j_t$ of $E|f_{t_j}(V)$ for
$t_j\le t\le t_{j+1}$  satisfying
$$ s\circ \xi^j_t(f_{t_j}(x)) = f_t(x) \quad (x\in V,\ t_j\le t\le t_{j+1}). $$
In particular we have
$s\circ \xi^{j}_{t_{j+1}}(f_{t_j}(x))=f_{t_{j+1}}(x)$ for $j=0,1,\ldots,k-1$.
Comparing these compatibility conditions with those defining the iterated
spray bundle $E^{(k)}$ (3.4) we see that these $k$ families can be joined into
a single family of sections $\xi_t$ ($0\le t\le 1$) of $E^{(k)}|f_0(V)$.
Explicitly we define for $x\in V$ and $t_j\le t\le t_{j+1}$:
$$
    \xi_t(f_0(x)) = \bigl( \xi^{0}_{t_1}(f_0(x)),\xi^{1}_{t_2}(f_{t_1}(x)),
               \ldots, \xi^{j-1}_{t_j}(f_{t_{j-1}}(x)),
               \xi^{j}_{t}(f_{t_j}(x)),0,\ldots,0\bigr) \in E^{(k)}
$$
(the last $k-j-1$ components are the zero elements in the fiber of $E$
over $s\circ \xi^{j}_{t}(f_{t_j}(x))=f_t(x)$). One easily verifies that
these sections satisfy the stated conditions.
\endpr

%
%
%
%
\beginsection 4. The h-Runge theorems.

In this section we prove Runge-type approximation theorems for
\holo\ sections of submersions with a spray over a Stein base.
Theorems 4.1 and 4.2 are due to Gromov [Gro], while theorem 4.5 is new.
The idea in these proofs is essentially the same as in Grauert's paper
[Gr2] where the analogous results had been proved for sections of
bundles over Stein spaces whose fiber is a complex Lie group $L$ and
the transition functions are Lie group automorphisms of $L$.

To motivate the discussion we recall that, if $X$ is a Stein manifold and
$K\subset X$ is a compact \hc\ subset, we can approximate each function
\holo\ in a \nbd\ of $K$ uniformly on $K$ by functions holomorphic on $X$.
This is the Oka-Weil theorem [H\"o2] which extends the classical Runge theorem for
planar sets $K \subset \C$ with connected complement.
Of course we cannot expect such results for sections of an arbitrary holomorphic
submersion over $X$; in particular this fails for maps of Stein manifolds into a
hyperbolic complex manifold $F$. On the other hand, if $h\colon Z\to X$
admits a fiber-spray (def.\ 3.1) then the Runge approximation property
is homotopy independent in the following sense:
\it If $f_t$ is a homotopy of sections in a \nbd\ of $K$ such that $f_0$
has a \holo\ extension to $X$, then each section in the homotopy can
be approximated on $K$ by sections \holo\ on $X$. \rm
Following Gromov we call such results {\it h-Runge theorems}.
We first state the special case, theorem 4.1.
The general parametric case is explained in theorem 4.2.

%
%
%
%
\proclaim 4.1 Theorem:
Let $X$ be a Stein manifold and $h\colon Z\to X$ a holomorphic
submersion which admits a fiber-spray (def.\ 3.1).
Let $K\subset X$ be a compact \hc\ set. Assume that $U\subset X$
is an open set containing $K$ and $f_t\colon U \to Z$ ($0\le t\le 1$)
is a homotopy of holomorphic sections of $h\colon Z\to X$ over $U$
such that $f_0$ extends to a holomorphic section over $X$.
Then for any metric $d$ on $Z$ and any $\e>0$ there exists a
continuous family of holomorphic sections $\tilde f_t\colon X\to Z$
($0\le t\le 1$) such that $\tilde f_0=f_0$ and
$$
   d\bigl(\tilde f_t(x),f_t(x)\bigr) <
      \e \quad (x\in K,\ 0\le t\le 1).
$$

\ni\it Remark. \rm
In fact we will prove the following stronger result:
\sl There is a continuous family of continuous sections
$g_{t,u} \colon X\to Z$\  ($0\le t,u\le 1$)
which are holomorphic in a neighborhood of $K$ and satisfy:
\item{(a)} $g_{t,0}=f_t$ on $K$ for all $0\le t\le 1$,
\item{(b)} the section $\tilde f_t=g_{t,1}$ is holomorphic on $X$ for each $t\in [0,1]$,
\item{(c)} $g_{0,u}=f_0$ for all $0\le u\le 1$, and
\item{(d)} $d\bigl(g_{t,u}(x),f_t(x)\bigr) < \e$ for $x\in K$ and $0\le t,u\le 1$.

\ni\rm The existence of such a homotopy $g_{t,u}$ connecting $f_t$ and $\tilde f_t$
comes from the proof of theorem 4.1 and will be useful to us in sect.\ 6.

\demo Proof: Let $(E,p,s)$ be the spray on $Z$ associated to $h$.
After shrinking $U$ around $K$ we obtain by proposition 3.6 an integer $k>0$
and a homotopy of sections $\xi_t$ over the set $f_0(U)\subset Z$ of the iterated
spray bundle $p^{(k)}\colon E^{(k)} \to Z$ which are mapped back to $f_t$ by
$s^{(k)} \colon E^{(k)}\to Z$. In particular $\xi_0$ is the zero section.
Since $X$ is Stein, the restriction $E^{(k)}|f_0(X)$ admits the structure of
a \hvb\ by corollary 3.5. It now suffices to approximate $\xi_t$ on the \hc\
subset $f_0(K)$ of $f_0(X)$ by a homotopy of \holo\ sections of
$E^{(k)}|f_0(X)$ (keeping the zero section fixed) and to take $\tilde f_t$ to be
their images in $Z$ by the spray map $s^{(k)}$. This can be done by the usual
Oka-Weil approximation theorem for sections of a holomorphic vector bundle,
and it can be reduced to the aproximation of functions by embedding
the given bundle as a subbundle of a trivial bundle. Even though this is
standard, we outline the proof for later purposes.

For convenience we let $p'\colon E'=f_0^*(E^{(k)}) \to X$ be the
pull-back of $E^{(k)}$ by the section $f_0 \colon X\to Z$, and we denote by
$s'\colon E'\to Z$ the \holo\ map induced by the spray
$s^{(k)} \colon E^{(k)} \to Z$.
We may then consider $\xi_t$ as sections of the bundle $p'\colon E'\to X$
over $U\subset X$ such that $\xi_0$ is the zero section and
$s'\circ \xi_t(x)=f_t(x)$ for all $x\in U$ and $0\le t\le 1$.
Choose a smooth function $\chi\colon X\to [0,1]$ which is identically
one in a \nbd\ of $K$ and has compact support contained in $U$.
Since  $K$ is \hc\ in $X$, there is a smooth \psh\ exhaustion function
$\rho\colon X\to\R_+$ which vanishes in a \nbd\ $U_0 \ss U$ of $K$
and is strictly positive on the support of $d\chi$. We can choose
$\rho$ such that for each fixed value of $\tau>1$ there is a
section $v_t$ of $E'$ which solves the equation
$$
    \dibar v_t = \dibar(\chi\xi_t)= \xi_t\, \dibar \chi            \eqno(4.1)
$$
and whose $L^2$ norm with weight $e^{-\tau\rho}$ (measured in a fixed hermitian
metric) is bounded on each compact set in $X$  by a constant times the norm of
the data $\xi_t \dibar\chi$ (with a constant independent of $\tau$).
By H\"ormander [H\"o1] such a solution is given by a linear operator
$v_t= T_\tau(\xi_t\dibar \chi)$. We have $v_0=0$, each $v_t$ is smooth,
and the family is continuous in $t$. Set
$$
   \tilde g_{t,u} = \chi \xi_t - uv_t,\quad g_{t,u}=s'\circ \tilde g_{t,u}
   \quad (0\le t,u\le 1).
$$
Clearly $g_{t,u} \colon X\to Z$ is a continuous family of sections of $Z\to X$ which
satisfies (a)--(c) in the remark following theorem 4.1. By choosing $\tau$
sufficiently large (depending on $\e$) the family $g_{t,u}$ will also satisfy (d)
which can be seen as follows. When $\tau\to\infty$, the $L^2$ norm with
weight $e^{-\tau\rho}$ of $\xi_t\dibar \chi$ tends to zero since $\rho>0$ on
$\supp \dibar\chi$. Since $\rho$ vanishes in $U_0$, it follows that
unweighted $L^2(U_0)$ norm $||v_t||_{L^2(U_0)}$ tends to zero as $\tau\to+\infty$.
By the Cauchy estimates the sup norm of $v_t|K$ tends to zero and hence the
sections $g_{t,u}$ converge to $f_t$ as $\tau\to +\infty$, uniformly in $t$ and $u$.
In particular, $\tilde f_t=g_{t,1}$ ($0\le t\le 1$) is a homotopy of holomorphic
sections of $Z\to X$ satisfying theorem 4.1.
\endpr

We will also need the following parametric version of the h-Runge theorem.
%
%
%
%
\proclaim 4.2 Theorem:  Let $X$ be a Stein manifold and $h\colon Z\to X$ a
holomorphic submersion which admits a spray (def.\ 3.1). Let $K\subset X$
be a compact \hc\ set and let $U, V\subset X$ be open, relatively compact
subsets in $X$ such that $K\subset U\subset V$.
Let $Y$ be a compact Hausdorff space (the parameter space),
$Y_0 \subset Y$ a compact subset, and $Y'\subset Y$ an open set containing $Y_0$.
Assume that $f_{y,t}\colon U\to Z$ is a family of \holo\ sections
of $h\colon Z\to X$, depending continuously on $y\in Y$ and $0\le t\le 1$,
such that the sections $f_{y,0}$ ($y\in Y$) and $f_{y,t}$ ($y\in Y'$, $0\le t\le 1$)
extend to holomorphic sections over $X$. Then for each $\e>0$ there exists a
continuous family of holomorphic sections $\tilde f_{y,t}\colon V\to Z$ over $V$
($y\in Y,\ 0\le t\le 1$) satisfying
\item{(a)} $\tilde f_{y,0}=f_{y,0}$ for all $y\in Y$,
\item{(b)} $\tilde f_{y,t} = f_{y,t}$ for all $y\in Y_0$ and $0\le t\le 1$, and
\item{(c)} $d\bigl( \tilde f_{y,t}(x), f_{y,t}(x) \bigr) <\e$ for all
$x\in K$, $y\in Y$ and $0\le t\le 1$.

\demo Remarks:
1. As in theorem 4.1 the proof will show that the approximating family
$\tilde f_{y,t}$ can be chosen so that it can be connected to the initial
family $f_{y,t}$ by a homotopy of sections $g_{y,t,u}$ $(u\in [0,1])$ such that
the homotopy is fixed for each $y\in Y_0$ (where $\tilde f_{y,t}=f_{y,t}$).

\ni 2. By a standard limiting procedure we could obtain a
family $\tilde f_{y,t}$ on all of $X$, but we shall not need this.

\demo Proof:
We begin by reducing to the approximation problem for families of
sections of an iterated spray bundle. This is essentially proposition 3.6 with
the addition of the parameter $y\in Y$. Let $(E,p,s)$ be a spray on $Z$ and
$(E^{(k)},p^{(k)},s^{(k)})$ its $k$-th iterated spray (3.4).

\medskip\ni\bf 4.3 Proposition. \rm (Assumptions as in theorem 4.2.) \sl
Let $U'$ be an open set in $X$ such that $K\subset U'\ss U$.
Then there are an integer $k>0$ and a continuous family of \holo\ sections
$\xi_{y,t}$ of $E^{(k)}|f_{y,0}(U')$\ \ ($y\in Y,\ 0\le t\le 1$) such that
\item{(i)}    $\xi_{y,0}$ is the zero section for each $y\in Y$,
\item{(ii)}   $\xi_{y,t}$ extends to a \holo\ section of $E^{(k)}|f_{y,0}(V)$
for each $y$ in a \nbd\ of $Y_0$ and $t\in [0,1]$, and
\item{(iii)}  $s^{(k)} \circ \xi_{y,t}(f_{y,0}(x)) = f_{y,t}(x)$ for
    all $x\in U'$, $y\in Y$ and $0\le t\le 1$. Moreover, for each $y$ in a
    \nbd\ of $Y_0$ this holds for all $x\in V$.

\rm
\demo Proof: It suffices to prove that each fixed $t_0\in [0,1]$ has a
\nbd\ $I_0 \subset [0,1]$ such that there exists a family of \holo\ sections
$\xi_{y,t}$ of the vector bundle $E|f_{y,t_0}(U')$, depending continuously
on $y\in Y$ and $t\in I_0$, such that $\xi_{y,t_0}$ is the zero section and
$$ s\circ \xi_{y,t}(z) =f_{y,t}(h(z)) \quad
                   (z\in f_{y,t_0}(U'),\ y\in Y,\ t\in I_0).               \eqno(4.2)
$$
Moreover, for $y$ in a \nbd\ of $Y_0$ the property (4.2) must hold over the larger set
$z\in f_{y,t_0}(V)$. Proposition 4.3 then follows from this as in the proof of
proposition 3.6 by using the compactness of $[0,1]$ and combining the finitely
many families of sections $\xi^j_{y,t}$ obtained over $k$ subintervals
$t\in [t_j,t_{j+1}] \subset [0,1]$ into a single family of sections of the
iterated spray bundle $E^{(k)}$ over the sets $f_{y,0}(U')$ resp.\ $f_{y,0}(V)$.

We will consider the case $t_0=0$ when all initial sections $f_{y,0}$ exist
over $V$. The only difference for $t_0>0$ is that some sections only exist
over the smaller set $U'$, but the proof goes through in the same way.
Denote by $E_0=\ker VDs \subset E$ the kernel of the vertical derivative of $s$ (3.2)
and let $\pi \colon E\mapsto \tilde E = E/E_0$ be the quotient projection.
Recall that a {\it holomorphic splitting} of $\pi$ is a \hvb\ homomorphism
$G\colon \tilde E\to E$ such that $\pi\circ G$ is the identity on $\tilde E$;
in such case we have $E=E_0\oplus G(\tilde E)$. A splitting exists over any
Stein subset of $Z$ [GRo, p.\ 256].

\proclaim 4.4 Lemma: {\rm (Hypotheses as in theorem 4.2.)}
There is a family of  holomorphic splittings
$G_y\colon \tilde E|f_{y,0}(V) \to E|f_{y,0}$
which depends continuously on $y\in Y$, and hence there is a holomorphic
direct sum splitting
$$ E|f_{y,0}(V)= E_0|f_{y,0}(V) \oplus E'_y                             \eqno(4.3) $$
depending continuously on $y\in Y$.

\demo Remark: We may consider the restricted bundles in lemma 4.4 as subsets of
the bundle $E$ resp.\ $\tilde E=E/E_0$, and the continuity of the family $G_y$
with respect to $y$ should be understood in this sense. Lemma 4.4 also holds
if we replace $t=0$ by an arbitrary $t_0\in [0,1]$, except that in this case
one must replace for each $y\in Y\bs Y'$ the set $V$ in (4.3) by the smaller
set $U'$ (since the section $f_{y,t_0}$ is only defined on $U'$).

\demo Proof of lemma 4.4: For each fixed $y\in Y$ the section $f_{y,0}(X) \subset Z$
(which is a Stein submanifold of $Z$) is contained in an open Stein set
$D_y\subset Z$ according to a theorem of Siu [Siu]. Hence by compactness of $Y$
there is an open covering $\{Y_j \colon 1\le j\le J\}$ of $Y$ and a family of
open Stein subsets $\{D_j\colon 1\le j\le J\}$ of $Z$ such that
$f_{y,0}(V)\subset D_j$ when $y\in Y_j$. Let $\chi_j$ be a continuous
partition of unity on $Y$ subordinate to the covering $\{Y_j\}$.
Let $H_j\colon \tilde E|D_j \to E|D_j$ be a splitting of $\pi\colon E\to \tilde E$
over $D_j$ (such $H_j$ exists since $D_j$ is Stein). Then the family
$$ G_y = \sum_{j=1}^J \chi_j(y) H_j \colon \ \
             \tilde E|f_{y,0}(V) \to E|f_{y,0}(V) $$
satisfies lemma 4.4. Note that the map $G_y$ is well defined since the
coefficient $\chi_j(y)$ vanishes when $f_{y,0}(V)$ is not contained in $D_j$,
and $G_y$ is a splitting of $\pi$ since it is a convex
linear combination of finitely many splittings.
\endpr

We continue with the proof of proposition 4.3. By the inverse function theorem
the map $s\colon  E'_y \to Z$ (the restriction of the spray map $s\colon E\to Z$
to $E'_y$) is a biholomorphic map from a \nbd\ of the zero section of $E'_y$
onto a \nbd\ of the set $f_{y,0}(V)$ in $Z$ (since the vertical derivative $VDs$
is an isomorphism when restricted to $E'_y$). Denote the local inverse of this map
by $u_y$. Because of continuous dependence on the data on $y\in Y$ and the compactness
of $Y$ the neighborhoods on which the maps $u_y$ are defined can be chosen uniformly
with respect to $y\in Y$. Hence there is a $t_1>0$ such that for $0\le t\le t_1$
and $y\in Y$  the set $f_{y,t}(U')$ belongs to the domain of $u_y$. The sections
$\xi_{y,t}=u_y\circ f_{y,t}\circ h \colon f_{y,0}(U') \to E'_y \subset E$
for $0\le t\le t_1$ then satisfy (4.2). This completes the proof
of proposition 4.3.
\endpr

To complete the proof of theorem 4.2 it remains to approximate
the sections $\xi_{y,t}$ uniformly on $f_{y,0}(K)$ by \holo\ sections
of $E^{(k)}|f_{y,0}(V)$. We shall first do this locally with respect to
$y$ and finally patch the approximations together. Let $Y'\supset Y_0$ be a \nbd\
of $Y_0$ such that the sections $\xi_{y,t}$ for $y\in Y'$ are defined and
\holo\ over the set $f_{y,0}(V) \subset Z$. As in the proof of lemma 4.4
there exist a finite open covering $\{Y_j \colon 1\le j\le J\}$ of $Y$
and open Stein subsets $D_j \subset Z$ such that $f_{y,0}(V) \subset D_j$ for
all $y$ in a \nbd\ $Y'_j\supset Y_j$ of $\bar Y_j$. By corollary 3.5 the
restriction $E^{(k)}|D_j$ has a vector bundle structure and we can embed it as
a vector subbundle of a trivial bundle $D_j\times \C^N$. With this identification
we may write
$$ \xi_{y,t}( f_{y,0}(x)) = \bigl( f_{y,0}(x), g_{y,t}(x) \bigr)  $$
where $g_{y,t}\colon U'\to \C^N$ is a \holo\ map. Note that
$g_{y,0}=0$ for all $y\in \bar Y_j$, and $g_{y,t}$ extends  holomorphically to $V$
for $y\in Y'\cap Y'_j$. Choose a continuous function $u\colon Y\to [0,1]$
which vanishes in smaller \nbd\ of $Y_0$ and is identically one outside $Y'$.
In now suffices to approximate the functions $g_{y,t}$ by functions \holo\
in $V$ without changing their values for $y\in Y'$ (when they are already
holomorphic in $V$). This is done by taking
$$ \tilde g_{y,t}(x) = u(y)\chi(x) g_{y,t}(x) - v_{y,t}(x), $$
where $v_{y,t}$ is the solution of the $\dibar$-equation
$$
   \dibar_x v_{y,t}(x) = u(y) g_{y,t}(x) \dibar_x \chi(x)
$$
with minimal $L^2$-norm with respect to a suitable weight and
$\chi\colon X\to\R_+$ is a cut-off function as in the proof of
theorem 4.1. Since $u(y)=0$ for $y$ near $Y_0$,
we have $v_{y,t}=0$ and hence $\tilde g_{y,t}=g_{y,t}$ for such $y$.

This gives for each $j=1,2,\ldots,J$ an approximating family of \holo\ sections
$\xi^j_{y,t}$ of $E^{(k)}|f_{y,0}(V)$ for $y$ in a \nbd\ $Y'_j$ of $\bar Y_j$.
To conclude the proof of theorem 4.2 it suffices to combine these families  into
a single family of sections $\tilde \xi_{y,t}$ ($y\in Y$) and then take
$$ \tilde f_{y,t}=s^{(k)} \circ \tilde \xi_{y,t} \circ f_{y,0}
      \qquad (y\in Y,\ 0\le t\le 1).                          \eqno(4.4)$$
However we cannot do the patching by the usual partition of unity in the
$y$ variable since the vector bundle structures on different restrictions
$E^{(k)}|D_j$ may not agree on their intersection. Instead we can do a
stepwise extension as follows. For $y\in Y'_1$ we take
$\tilde\xi_{y,t}=\xi^1_{y,t}$. Choose a continuous function
$\chi_1 \colon Y \to [0,1]$ such that $\chi_1(y)=1$ for $y\in \bar Y_1$ and
$\supp\chi_1 \subset Y'_1$. For $y\in Y'_2$ we then set
$$ \tilde \xi_{y,t}= \chi_1(y) \xi^1_{y,t} + (1-\chi_1(y)) \xi^2_{y,t}  \eqno(4.5) $$
where the linear combination is taken with respect to a vector bundle
structure on $E^{(k)}|D_2$. This definition is good since $\chi_1(y)=0$
for those values of $y$ (i.e., for $y \in Y\bs Y'_1$) for which the section
$\xi^1_{y,t}$ is not defined, and so the patching only occurs over the set
$y\in Y'_2$. Clearly the family (4.5) extends continuously to all parameter
values $y\in Y'_1\cup Y'_2$ and it equals the previously chosen family
of sections for $y\in Y_1$. Moreover, for $y\in Y'_1\cap Y'_2$ both sections
$\xi^1_{y,t}$ and $\xi^2_{y,t}$ approximate the initial section $\xi_{y,t}$
over $f_{y,0}(K)$, and hence the same is true for their convex combination (4.5).
We now continue in the same way by patching the family (4.5) with $\xi^3_{y,t}$
over the parameter set $y\in Y'_3$ with respect to the vector bundle structure
on $E^{(k)}|D_3$. In finite number in steps we obtain a continuous
family of sections $\tilde \xi_{y,t}$ such that the sections (4.4) satisfy
theorem 4.2.
\endpr

Our next result, which is an immediate application of theorem 4.2,
is essential in our approach to extending a holomorphic section
across a \pb\ in the critical case.

\proclaim 4.5 Theorem:  Let $X$ be a Stein manifold and
$Z=X\times F\to X$ a trivial bundle whose fiber $F$
admits a spray (def.\ 1.2).
Suppose that $U\subset X$ is a Stein domain which is Runge
in $X$ and $\Theta_u\colon X\to X$ is a family of holomorphic
mappings, depending continuously on $u\in [0,1]$, such that
\item{(i)}   $\Theta_0$ is the identity on $X$,
\item{(ii)}  $\Theta_u(U)\subset U$ for all $u\in [0,1]$, and
\item{(iii)} $\Theta_1(X)\subset U$.

\ni\sl Let $f_t\colon U\to Z$ ($0\le t\le 1$) be a homotopy of
holomorphic sections such that $f_0$ and $f_1$ extend holomorphically to
$X$. Let $d$ be any metric on $Z$. For each choice of a compact set
$K\subset U$, a relatively compact set $V\ss X$, and $\e>0$ there exists
a homotopy of holomorphic sections
$\tilde f_t\colon V \to Z$ ($0\le t\le 1$) such that
$\tilde f_0=f_0$, $\tilde f_1=f_1$, and
$$
    d\bigl( \tilde f_t(x), f_t(x) \bigr) <
       \e \quad (x\in K,\ 0\le t\le 1).
$$
\rm

\demo Remarks: 1. In a typical application of theorem 4.5 the sets
$U\subset X$ are bounded convex domains in $\cn$ and $\Theta_u$ is a
family of linear contractions to a point in $U$. By a limiting argument
it is possible to prove that such a holomorphic homotopy $\tilde f_t$
exists on all of $X$, but we shall not need this.

\ni 2. The analogous result holds, with the same proof, for
parametrized families of sections. We leave out the obvious
details.

\demo Proof: By reparametrizing the family $f_t$ we may assume
that for some small $\delta>0$ we have $f_t=f_0$ for $0\le t\le \d$
and $f_t=f_1$ for $1-\d \le t\le 1$. Choose a continuous function
$u\colon [0,1]\to [0,1]$ such that $u(t)=0$ for $t$ near $0$ or $1$,
and $u(t)=1$ for $\d\le t\le 1-\d$. We shall identify sections
of the trivial bundle $X\times F \to X$ with mappings $X\to F$. Set
$$ f_{t,s}= f_t \circ \Theta_{(1-s)u(t)} \quad (0\le t,s\le 1). $$
This family satisfies the following properties:

\item{(a)} $f_{t,s}$ is defined and holomorphic in $U$ for each
$0\le t,s\le 1$,

\item{(b)} $f_{t,0}=f_t \circ \Theta_{u(t)}$ is defined and holomorphic
on all of $X$ for each $t\in [0,1]$,

\item{(c)} $f_{t,1}=f_t\circ \Theta_0 = f_t$ for each $t\in [0,1]$, and

\item{(d)} $f_{0,s}=f_0\circ \Theta_{u(0)} = f_0$ and
           $f_{1,s}=f_1\circ \Theta_{u(1)} = f_1$ for all $s\in [0,1]$.

\ni
It remains to apply theorem 4.2 with the parameter space $t\in Y=[0,1]$ and
the subspace $Y_0=\{0,1\}$. (Note that our current variable $s$ plays the role
of the time parameter $t$ in theorem 4.2.) If $\tilde f_{t,s}$ is the approximating
family of holomorphic sections in $V$ as in theorem 4.2 then the sections
$\tilde f_t=\tilde f_{t,1}$ ($0\le t\le 1$) satisfy theorem 4.5.
\endpr

%
%
%
%
\beginsection 5. Gluing holomorphic sections over Cartan pairs.

The main results of this section are theorems 5.1 and 5.5 on
gluing holomorphic sections over Cartan pairs. This can be
found in sect.\ 1.6 of [Gro] and is similar to the
corresponding results of Grauert [Gr3] and Cartan [Car].

%
%
%
%
\proclaim 5.1 Theorem:
Let $h\colon Z\to X$ be a holomorphic submersion onto a Stein
manifold $X$. Let $d$ be a metric on $Z$.
Let $(A,B)$ be a Cartan pair in $X$ (def.\ 2.3) such that
the set $C=A\cap B$ is Runge in $B$. Suppose that $\tilde B$ is an open
\nbd\ of $B$ in $X$ such that the restriction $Z|\tilde B = h^{-1}(\tilde B)$
admits a fiber-spray (def.\ 3.1). Let $\tilde A\subset X$ be an
open \nbd\ of $A$ and $a\colon \tilde A\to Z$ a \holo\ section of $Z\to X$
over $\tilde A$. Then for each $\e>0$ there is a $\d>0$ satisfying the
following property. If $b\colon \tilde B\to Z$ is a \holo\ section
satisfying $d\bigl( a(x), b(x)\bigr) <\d$ for $x\in \tilde C=\tilde A\cap \tilde B$,
there exist homotopies $a_t$ (resp.\ $b_t$), $0\le t\le 1$,
of holomorphic sections over a \nbd\ $A'$ of $A$ (resp.\ over a \nbd\ $B'$ of $B$)
such that $a_0=a$, $b_0=b$, $a_1=b_1$ on $C'=A'\cap B'$, and
$$ \eqalign{ d\bigl( a_t(x), a(x) \bigr) &< \e \quad (x\in A',\ 0\le t\le 1); \cr
    d\bigl( b_t(x), b(x) \bigr) &< \e \quad (x\in C',\ 0\le t\le 1). \cr}
$$

Using sprays we shall reduce the proof of theorem 5.1 to the model
case described by the following proposition; this is analogous to the
classical Cartan--Grauert attaching lemma [Gr3], [Car]. We denote by
$H^\infty(\Omega,\cn)$ the Banach space of bounded \holo\ maps
$\Omega\to\cn$ equipped with the sup norm over all components.

%
%
\proclaim 5.2 Proposition:
Let $(A,B)$ be a Cartan pair in a Stein manifold $X$ such that $C=A\cap B$ is
Runge in $B$. Let $\tilde C\subset X$ be an open \nbd\ of $C$, $U\subset \cn$ an
open \nbd\ of the origin in $\cn$, and $\psi_0\colon \tilde C\times U\to\cn$
a bounded \holo\ map such that for each $x\in \tilde C$, $\psi_0(x,0)=0$
and $\psi_0(x,\cdotp) \colon U\to\C^n$ is injective (i.e., biholomorphic onto its
image). Then there are \nbd s $A'\supset A$ and $B'\supset B$ with
$C'=A'\cap B'\ss \tilde C$, a \nbd\ $W$ of $\psi_0$ in the Banach space
$H^\infty(\tilde C\times U, \C^n)$,
and smooth Banach space operators ${\cal A'}\colon W\to H^\infty(A',\C^n)$,
${\cal B'}\colon W\to H^\infty(B',\C^n)$, with ${\cal A'}(\psi_0)=0$
and ${\cal B'}(\psi_0)=0$, such that for each $\psi\in W$ the
bounded \holo\ maps $\alpha={\cal A'}(\psi) \colon A'\to\cn$,
$\beta={\cal B'}(\psi) \colon B'\to \cn$ satisfy
$$ \psi\bigl(x, \alpha(x) \bigr)= \beta(x) \quad (x\in A'\cap B').      \eqno(5.1)$$
Moreover, if $\psi \in W$ satisfies $\psi(x,0)=0$ for $x\in\tilde C$
then ${\cal A'}(\psi)=0$ and ${\cal B'}(\psi)=0$.

\demo Remark: We can view a pair of maps satisfying (5.1) as a section of
a nonlinear bundle over $A'\cup B'$ obtained by patching the trivial bundles
over $A'$ resp.\ $B'$ by the map $\psi$. For later application to parametrized
families it is convenient to have a canonically given solution (i.e., by
operators), although this could be avoided by a suitable analogue
of Satz 8 in [Gr2].

\demo Proof of proposition 5.2:
By shrinking $\tilde C$ we may assume
that it is Runge in a \nbd\ $B_0$ of $B$.
We choose \nbd s $A'\supset A$ and $B'\supset B$
as in lemma 2.4 so that $B'\subset B_0$, $C'=A'\cap B'\ss \tilde C$,
and there are bounded linear operators
${\cal A}\colon H^\infty(C',\cn) \to H^\infty(A',\cn)$ and
${\cal B}\colon H^\infty(C',\cn) \to H^\infty(B',\cn)$ satisfying
$c= {\cal A}c - {\cal B}c$ for all $c\in H^\infty(C',\cn)$.
Consider first the case when $\psi_0(x,u)=u$ is the identity map in
the $u$-variable for each $x\in \tilde C$. Consider the operator
$$ \eqalign{ & \Phi\colon H^\infty(C',\cn)\times  H^\infty(\tilde C\times U, \C^n)
  \to H^\infty(C',\cn), \cr
    & \Phi(c,\psi)(x) = \psi\bigl( x,{\cal A}c(x) \bigr) - {\cal B}c(x)
    \qquad (x\in C').  \cr}
$$
We claim that $\Phi$ is defined and smooth for $c\in H^\infty(C',\cn)$ in a \nbd\ of
the origin and for $\psi\in H^\infty(\tilde C\times U, \C^n)$. Clearly $\Phi$ is linear
and hence smooth in $\psi$. To see that $\Phi$ is smooth in $c$ we choose a
\nbd\ $U'\subset \cn$ of $0$ such that $\bar U'\subset U$. By Cauchy estimates the
restriction map $\psi \to \psi|C'\times U'$ is a bounded linear operator from
the space $H^\infty(\tilde C\times U, \C^n)$ to $\cC^\infty(C'\times U',\cn)$.
On the set of $c$'s for which ${\cal A}c(x)\in U'$ for all $x\in C'$
(these form an open \nbd\ of the origin in $H^\infty(C',\cn)$) the first term in
$\Phi$ is the composition operator of a linear operator ${\cal A}$ with a
smooth map $\psi$. Hence $\Phi$ is a smooth operator. In fact we only need
that $\Phi$ is of class $\cC^1$ which is seen directly from the formula
for its differential
$$
    D\Phi(c,\psi)(c',\psi') = \psi'(\cdotp,{\cal A}c)+
                        D_2\psi(\cdotp,{\cal A}c) {\cal A}c' - {\cal B}c'.
$$
Note that $\Phi(c,\psi_0)= {\cal A}c - {\cal B}c= c$. Hence $D_c \Phi(0,\psi_0)$
(the partial derivative of $\Phi$ with respect to the first variable)
is the identity map on $H^\infty(C',\cn)$. By the implicit function theorem in
Banach spaces there is an open set $W\subset H^\infty(\tilde C\times U, \C^n)$
containing $\psi_0$ and a smooth map ${\cal C} \colon W\to H^\infty(C',\cn)$
such that $\Phi({\cal C}(\psi),\psi)=0$ and ${\cal C}(\psi_0)=0$. Moreover, if
$\psi\in W$ satisfies $\psi(x,0)=0$ for all $x\in \tilde C$, then
$c=0$ solves the equation $\Phi(c,\psi)=0$ and hence by local uniqueness of
solutions we have ${\cal C}(\psi)=0$ for any such $\psi$. The operators
${\cal A'} = {\cal A}\circ{\cal C}$ and ${\cal B'}={\cal B}\circ {\cal C}$
then satisfy proposition 5.2.

The general case (when $\psi_0$ is not the identity, or even close to the identity)
can be reduced to the special case as follows. Since $C$ is Runge in $B$,
there are open sets $\tilde C_0\subset X$, $B_0\subset X$, $U_0,U_1\subset \cn$,
satisfying $C\subset \tilde C_0\ss \tilde C$,
$0\in U_0\ss U_1\ss U$, $B\subset B_0$, such that we can
approximate $\psi_0$ as well as desired on $\tilde C_0\times U_1$ by a \holo\
map $\tilde \psi\colon B_0\times \cn\to \cn$ such that $\tilde \psi(x,\cdotp)$
is biholomorphic on $U_1$ for each $x\in \tilde C_0$ and $\tilde \psi(x,0)=0$
for all $x\in B_0$. If the approximation of $\psi_0$ by $\tilde \psi$ is
sufficiently close on $\tilde C_0\times U_1$, there is a unique \holo\ map
$\psi'_0\colon \tilde C_0\times U_0\to U_1$ which satisfies
$$ \psi_0(x,u)= \tilde\psi(x, \psi'_0(x,u))\quad (x\in \tilde C_0,\ u\in U_0) $$
and which is so close to the map $(x,u)\to u$ on the set $\tilde C_0\times U_0$
that it belongs to the domain $W\subset H^\infty(\tilde C_0\times U_0,\cn)$
of the operators ${\cal C}$, ${\cal A}'$, and  ${\cal B}'$  obtained in the special
case (with respect to the smaller set $\tilde C_0\times U_0$). We may assume that
the sets $A'\supset A$, $B'\supset B$, related to the operators
${\cal A}'$ resp.\ ${\cal B}'$ as above, satisfy $B'\subset B_0$
and $A'\cap B'\ss\tilde C_0$. Note that ${\cal C}(\psi'_0)=0$ since $\psi'_0(x,0)=0$.
Moreover, for each $\psi$ which is sufficiently uniformly close to $\psi_0$ on
$\tilde C\times U$ we have
$$ \psi(x,u)= \tilde\psi(x, \psi'(x,u)) \quad  (x\in \tilde C_0,\ u\in U_0) $$
where $\psi' \colon \tilde C_0\times U_0\to U_1$ belongs to $W$. In fact,
$\psi\to\psi'$ defines a smooth Banach space operator
$\Psi\colon W_0\subset H^\infty(\tilde C\times U,\cn)
    \to H^\infty(\tilde C_0\times U_0,\cn)$
in an open \nbd\ $W_0$ of $\psi_0$, with range in $W$. For each $\psi\in W_0$ and
$\psi'=\Psi(\psi) \in W$ the \holo\ maps $\a={\cal A}'(\psi') \colon A'\to\cn$
and $\b'={\cal B'}(\psi') \colon B'\to \cn$ satisfy $\psi'(x,\a(x))=\b'(x)$
for $x\in C'$. Hence the pair $\a(x)$ and $\b(x)=\tilde \psi(x,\b'(x))$
(the latter one is defined and \holo\ in $B'$) satisfies
$$ \psi(x,\a(x))= \tilde\psi\bigl( x,\psi'(x,\a(x))\bigr)
   = \tilde \psi(x,\b'(x)) = \b(x)
$$
which is precisely (5.1). By construction $\a$ and $\b$ are obtained from
$\psi$ by a composition of smooth Banach space operators.
\endpr

In order to reduce the proof of theorem 5.1 to proposition 5.2 we
also need the following lemma. Denote by $B^n(\eta)\subset \cn$ the open
ball of radius $\eta$.

%
%
\proclaim 5.3 Lemma: Let $h\colon Z\to X$ be a \holo\ submersion onto
a Stein manifold $X$, let $A\subset X$ be a compact set with a Stein \nbd\ basis,
and let $a\colon \tilde A\to Z$ be a \holo\ section defined in an open set
$\tilde A\supset A$. Then there are an integer $n>0$, a number $\eta>0$, a Stein
open set $U\subset Z$ containing $a(A)$, and a \holo\ map
$\tilde s\colon U\times B^n(\eta) \to Z$ such that for all $z\in U$ we have
\item{(i)}   $h\circ\tilde  s(z,t)=h(z)$ for all $t=(t_1,t_2,\ldots,t_n) \in B^n(\eta)$,
\item{(ii)}  $\tilde s(z,0)=z$, and
\item{(iii)} the vectors $\tilde V_j(z)= {\di\over \di t_j}\tilde s(z,0)$
($1\le j\le n$) span $VT_z(Z)$.

\demo Remark: Note that $\tilde s$ satisfies all requirements for a fiber-spray
except that it is not defined globally on $U\times \C^n$.

\demo Proof of Lemma 5.3:  The set $a(A)\subset Z$ has a basis of Stein \nbd s
in $Z$ according to [Siu] and [Shd]. By Cartan's theory there exist finitely many \hvf s
$\tilde V_1,\tilde V_2,\ldots,\tilde V_n$ in a \nbd\ of $a(A)$ in $Z$ which are tangent
to $VT(Z)$ (i.e., they are \holo\ sections of the vertical tangent bundle $VT(Z)$)
and which span $VT(Z)$ at each point in the given set. Let $\theta_j^t$ be
the flow of $\tilde V_j$. There is a small Stein \nbd\ $U \subset Z$ of $a(A)$ and
an $\eta>0$ such that the map $\tilde s\colon U\times B^n(\eta)\to\cn$, given by
$$
    \tilde s(z,t)= \tilde s(z,t_1,\ldots,t_n)=
    \theta_1^{t_1}\circ\theta_2^{t_2}\circ\cdots\circ\theta_n^{t_n}(z),
$$
satisfies all requirements.
\endpr

\demo Proof of theorem 5.1:   Let $\tilde s \colon U\times B^n(\eta) \to Z$ be the
local spray given by lemma 5.3. Choose a Stein open set $\tilde A\supset A$ so that
$a(\tilde A)\ss U$. Set $s_1(x,t) = \tilde s(a(x),t)$ for $x\in\tilde A$
and $t\in B^n(\eta)$. Then $s_1(x,0)=a(x)$ and $s_1$ is a submersion along $t=0$
(i.e., a local spray onto a \nbd\ of $a(\tilde A)$ in $Z$). Suppose that
$\tilde B\supset B$ is an open Stein set and $b\colon \tilde B\to Z$ is a
\holo\ section such that $b$ approximates $a$ in an open \nbd\ $\tilde C$
of $C=A\cap B$. Our goal is to construct a fiber preserving
\holo\ map $s_2\colon \tilde B\times \C^n\to Z$ such that $s_2(x,0)=b(x)$ and such that
$s_2$ is a submersion near $t=0$ which approximates $s_1$ in a \nbd\ of
$C\times \{0\}^n$. We then solve the equation
$s_2(x,\psi(x,t))=s_1(x,t)$ to get a map $\psi$ as in proposition 5.2.
If $\a$ and $\b$ satisfy proposition 5.2
(i.e., $\psi(x,\a(x))=\b(x)$ for $x$ near $C$) then
$$
        a_1(x)=s_1(x,\a(x)),\qquad  b_1(x)=s_2(x,\b(x)
$$
are \holo\ sections of $Z$ over \nbd s of $A$ resp.\ $B$ which agree near $C=A\cap B$.

To construct $s_2$ we assume that $\tilde B\supset B$ is so small that
there exists a fiber-spray $(E,p,s)$ over $Z|\tilde B=h^{-1}(\tilde B)$.
Choose a Stein open set $U' \subset Z$ such that
$a(C)\subset U'\ss U \cap h^{-1}(\tilde B)$.
Let $\tilde V_j$ be the \vf s as in lemma 5.3. Since $VDs(0)\colon E\to VT(Z)$ is a
surjective vector bundle homomorphism and the set $U'$ is Stein, we can split
$E|U'=\ker VDs(0)\oplus E'$ and lift the sections $\tilde V_j$ to \holo\ sections
$V_j$ of $E' \subset E|U'$ such that $VDs(0_z)V_j(z)=\tilde V_j(z)$ for all
$z\in U'$ and $1\le j\le n$.  For each $z\in U'$ and each collection of
$n$ vectors ${\cal W}=\{W_1,\ldots,W_n\} \subset E_z$ we define a map
$s_{\cal W}(z,\cdotp) \colon \C^n \to Z_{h(z)}$ by
$$ s_{\cal W}(z,t)= s\biggl( \sum_{j=1}^n t_j\, W_j \biggr).                $$
We have $s_{\cal W}(z,0)=z$ and
${\di\over\di t_j}s_{\cal W}(z,0)=VDs(0_z) W_j$.
In particular, for the collection ${\cal V}(z)= \{V_1(z),\ldots,V_n(z)\}$ we get
for $j=1,2,\ldots,n$ and $z\in U'$
$$
    {\di\over\di t_j}s_{{\cal V}(z)}(z,0)= VDs(0_z) V_j(z) =
     \tilde V_j(z) = {\di\over\di t_j} \tilde s(z,0).                 \eqno(5.2)
$$
The map $s_{\cal W}(z,t)$ is \holo\ in all arguments, including ${\cal W}$.

\proclaim 5.4 Lemma: Let $d$ be a metric on $Z$ and $d'$ a metric on $E$.
There are numbers $\eta>0$, $\d>0$ with the following property.
For each pair of points $z,w\in U'$ such that $h(z)=h(w)$ and $d(z,w)<\d$,
and for each collection ${\cal W}=\{W_1,\ldots,W_n\} \subset E_w$ such that
$d'(W_j,V_j(w))<\d$ for $j=1,\ldots,n$, there is an injective holomorphic map
$\phi_{\cal W}(z,w,\cdotp) \colon B^n(\eta) \to \C^n$
satisfying
\item{(i)}   $s_{\cal W}(w,\phi_{\cal W}(z,w,t))= \tilde s(z,t)$,
\item{(ii)}  $\phi$ is \holo\ in all arguments $z,w,{\cal W},t$, and
\item{(iii)} $\phi_{\cal W}(z,z,0)=0$.

\demo Proof:
Since $U\subset Z$ is Stein, we have a splitting $U\times \C^n= M\oplus N$ where
$M_z$ is the kernel of $D_t \tilde s(z,0)$ (the $t$-derivative of $\tilde s$ at the
zero section) and $N$ is some holomorphic complementary bundle. We split the fiber
vectors $t=(t',t'') \in M_z\oplus N_z$ accordingly (so the splitting depends on
the base point $z$). For each $z\in U$ the restriction of
$\tilde s$ to the fiber $N_z$ maps a \nbd\ of $0_z\in N_z$ biholomorphically
onto a \nbd\ of $z$ in the fiber $Z_{h(z)}$. The same is true for the restriction
$$  t''\in N_z \to  \tilde s(z,(t',t''))\in Z_{h(z)}                  \eqno(5.3) $$
of $\tilde s$ to fibers $\{t'\}\oplus N_z$ for all sufficiently small
vectors $t'\in M_z$.

Now (5.2) shows that for each pair of points $z,w\in U'$ in the same fiber
$Z_{h(z)}$ which are sufficiently close together, for each sufficiently small
vector $t'\in M_z$, and for each  collection of $n$ vectors
${\cal W}=\{W_1,\ldots,W_n\} \subset E_w$ which are sufficiently close to the
corresponding vectors ${\cal V}(w)= \{V_1(w),\ldots,V_n(w)\}$, the map
$$ t''\in N_z \to s_{\cal W}(w,(t',t''))  \in Z_{h(z)}                   \eqno(5.4) $$
takes a \nbd\ of ${0}^{''}_z\in N_z$ in $N_z$ biholomorphically onto a \nbd\ of $w$
in $Z_{h(z)}=Z_{h(w)}$ such that the image also contains the point $z$. For such choice
of points and vectors we take $\phi'_{\cal W}(z,w,t',\cdotp) \colon N_z \to N_z$ to
be the map (5.3) followed by the (unique!) local inverse of (5.4) at $t''=0$,
and then take
$$ \phi_{\cal W}(z,w,(t',t''))= \bigl(t', \phi'_{\cal W}(z,w,t',t'') \bigr). $$
This map is defined for $t=(t',t'') \in M_z\oplus N_z=\cn$ in some \nbd\ of
the origin $0\in \C^n$ which we may take to be independent of $z$, $w$,
$\cal W$, provided that all conditions regarding closeness are satisfied.
It is easily verified that this map satisfies all required properties.
Since both maps (5.3) and (5.4) depend \holo ally on all arguments, so does $\phi$.
\endpr

Suppose now that $\tilde C\subset \tilde A\cap \tilde B$ is an open Stein set
containing $C$ which is Runge in $\tilde B$. Assume that
$b\colon \tilde B\to Z$ is a \holo\ section such that $b(\tilde C)\subset U'$,
where $U'\subset Z$ is the \nbd\ of $a(C)$ chosen above. We consider the
restrictions $V_j|b(\tilde C)$ as \holo\ sections of the
bundle $E|b(\tilde B)$ over the set $b(\tilde C)$. By shrinking $\tilde C$ we get
holomorphic sections $W_j$ $(1\le j\le n)$ of $E|b(\tilde B)$ which approximate the
sections $V_j$ uniformly on $b(\tilde C)$ as close as desired. Write
${\cal W}(z)=\{W_1(z),\ldots,W_n(z)\}$ for $z=b(x)\in b(\tilde B)$.
The maps $s_1\colon \tilde A\times B^n(\eta) \to Z$ and
$s_2\colon \tilde B\times \C^n \to Z$, given by
$$ \eqalignno{ s_1(x,t) &= \tilde s(a(x),t),  \cr
               s_2(x,t) &= s_{{\cal W}(b(x))} (b(x),t),              &(5.5) \cr}
$$
are holomorphic. If $b$ is sufficiently uniformly close to $a$ over $\tilde C$ and
if $W_j(b(x))$ is sufficiently close to $V_j(b(x))$ for each $j=1,\ldots,n$
and $x\in \tilde C$, then by lemma 5.4 we have for each $x\in \tilde C$
an injective \holo\ map
$$
    \psi(x,\cdotp)= \phi_{{\cal W}(b(x))} (a(x),b(x),\cdotp) \colon B^n(\eta)\to \C^n
$$
which solves the equation
$$ s_2(x,\psi(x,t)) = s_1(x,t)\quad (x\in C',\ t\in B^n(\eta)).             $$
Moreover, if the approximations are sufficiently close,
$\psi$ is uniformly close to the map
$$
    \psi_0(x,t)= \phi_{{\cal V}(a(x))} (a(x),a(x),t)
$$
which satisfies $\psi_0(x,0)=0$ for $x\in \tilde C$. By proposition 5.3
we get open sets $A'$, $B'$, $C'=A'\cap B'$ in $X$, with
$A\subset A'\subset \tilde A$ and $B\subset B'\subset \tilde B$,
and \holo\ maps $\a\colon A'\to \C^n$, $\b\colon B'\to\C^n$, such that
$\psi(x,\a(x))=\b(x)$ for $x\in C'$. The homotopies of sections
$$ a_t(x)=s_1(x,t\a(x)) \ \ (x\in A'),\qquad b_t(x)=s_2(x,t\b(x)) \ \ (x\in B')
                                                                    \eqno(5.6) $$
for $0\le t\le 1$ then satisfy theorem 5.1.

Perhaps a word is in order regarding the proximity of the sections $a_t$ to $a=a_0$
and of $b_t$ to $b=b_0$. The rate of approximation of $\psi_0$ by $\psi$
depends on the proximity of $b(x)$ to $a(x)$ and on the proximity of the vector
fields ${\cal W}(b(x))$ to the fields ${\cal V}(a(x))$ for $x\in \tilde C$.
This in turn determines the estimates on the norms
$||\a||_{H^\infty(A')}$ and $||\b||_{H^\infty(B')}$ (by proposition 5.2).
Since the map $s_1$ (5.5) only depends on the section $a$, the definition
(5.6) shows that the estimate of $d\bigl(a_t(x),a(x)\bigr)$ for $x\in A'$ and
$t\in [0,1]$ depends only on $||\a||_{H^\infty(A')}$, and we get the stated
approximation result over $A'$.
However, the map $s_2$ (5.5) depends both on $b$ and on $\cal W$.
These quantities are under control only on $\tilde C$ and not on all of
$\tilde B$ (since we apply Runge approximation). Therefore we can estimate
$d\bigl(b_t(x),b(x)\bigr)$ in terms of $||\b||_{H^\infty(B')}$ only
for points $x\in C'$ and not on $B'\bs C'$.
\endpr

The following is an extension of theorem 5.1 to parametrized families of
sections.

%
%
%
%
\proclaim 5.5 Theorem:
Let $h\colon Z\to X$ be a holomorphic submersion onto a Stein manifold $X$,
let $d$ a metric on $Z$, and let $(A,B)$ be a Cartan pair in $X$ (def.\ 2.3) such that
the set $C=A\cap B$ is Runge in $B$. Suppose that $\tilde B$ is an open \nbd\
of $B$ in $X$ such that the restriction $Z|\tilde B = h^{-1}(\tilde B)$
admits a fiber-spray (def.\ 3.1). Let $Y$ be a compact Hausdorff space
(the parameter space) and $Y_0\subset Y$ a compact subset.
Let $\tilde A\subset X$ be an open \nbd\ of $A$ and
$a\colon \tilde A \times Y\to Z$ a continuous map such that
for each $y\in Y$, $a(\cdotp,y)$ is a \holo\ section of $Z$ over $\tilde A$.
Then for each $\e>0$ there is a $\d>0$ satisfying the following property.
If $b\colon \tilde B \times Y \to Z$ is a continuous map such that for each
$y\in Y$, $b(\cdotp,y)$ is a \holo\ section of $Z$ over $\tilde B$
satisfying
$$ \eqalign{ & d\bigl( a(x,y), b(x,y)\bigr) <\d \quad
    (x\in \tilde C=\tilde A\cap \tilde B,\ \ y\in Y), \cr
    & a(x,y)=b(x,y) \quad\quad\quad\ (x\in \tilde C,\ \ y\in Y_0), \cr}
$$
then there exist smaller \nbd s $A'\supset A$, $B'\supset B$ and
homotopies $a_t \colon A'\times Y\to Z$ resp.\ $b_t \colon B'\times Y\to Z$
($0\le t\le 1$) of families of holomorphic sections such that
$a_0=a$, $b_0=b$, $a_1=b_1$ on $C'=A'\cap B'$, and  for each $0\le t\le 1$ we have
$$ \eqalign{ & d\bigl( a_t(x,y), a(x,y) \bigr) < \e \quad
            (x\in A',\ y\in Y), \cr
    & d\bigl( b_t(x,y), b(x,y) \bigr) < \e \quad\, (x\in C',\ y\in Y), \cr
    & a_t(x,y) = b_t(x,y)\quad\quad\quad (x\in C',\ y\in Y_0). \cr}
$$

\demo Proof: This can done by essentially repeating the proof of theorem 5.1
with the addition of the parameter $y$. We shall only indicate a few critical
places in the proof where it is not completely obvious what must be done.
First one needs for each $y\in Y$ \holo\ maps
$$ s_{1,y}\colon \tilde A\times B^n(\eta) \to Z,\qquad
   s_{2,y}\colon \tilde B\times \C^n \to Z                       \eqno(5.7)
$$
for some integer $n>0$ and some $\eta>0$ which are related to the sections
$a_y=a(\cdotp,y)$ resp.\ $b_y=b(\cdotp,y)$ as in (5.5) and
which depend continuously on the parameter $y\in Y$. To get $s_{1,y}$
we cover $Y$ by finite number of open sets $Y_j$ ($1\le j\le J$) such that
for each $j$ there is an open Stein set $U_j\subset Z$, with
$a_y(\tilde A)\subset D_j$ for each $y\in Y_j$, and there are
finitely many \hvf s $\tilde V^j_k$ ($1\le k\le k_j$) which generate
$VT(Z)$ at each point of $U_j$. (This is similar to the proof of lemma 4.4 above.)
Let $\chi_j\colon Y\to [0,1]$ ($1\le j\le J$)
be a continuous partition of the unity on $Y$ subordinate to the covering $\{Y_j\}$.
Then the vector fields $\chi_j(y) \tilde V^j_k$ ($1\le j\le J$, $1\le k\le k_j$)
are well defined and \holo\ in a \nbd\ of the section $a_y(\tilde A)$ for each
$y\in Y$ (since $\chi_j(y)=0$ for those $y\in Y$ for which the field $\tilde V^j_k$
is not defined near $a_y(\tilde A)$, i.e., for $y$ outside $Y_j$). Together
these $n$ fields generate the vertical tangent bundle $VT(Z)$ at each point $a_y(x)$
for $x\in \tilde A$ and $y\in Y$. Using these fields and their
local flows we get as before a family of submersions $s_{1,y}$ (5.7)
depending continuously on $y\in Y$.

With a similar argument (see lemma 4.4) we obtain a family of maps
$s_{2,y}$ (5.7) which are submersions in a \nbd\ of the zero section
over the set $\tilde C$ and which approximate $s_{1,y}$ near $\tilde C\times \{0\}^n$.
This gives a family of transition mappings $\psi_y$ as in proposition 5.2
which are continuous in $y$ and approximate a certain initial family
$\psi_{0,y}$. By proposition 5.2 we obtain families of bounded \holo\ maps
$\a_y\colon A' \to \C^n$, $\b_y\colon B'\to\cn$, depending continuously on
$y\in Y$ and satisfying
$$ \psi_y\bigl(x, \a_y(x)\bigr) = \b_y(x)\quad (x\in C',\ y\in Y). $$
This gives for all $y\in Y$ homotopies of sections
$$ a_t(x,y)=s_{1,y}\bigl(x,t\a_y(x)\bigr) \ \ (x\in A'),\qquad
   b_t(x,y)=s_{2,y}\bigl(x,t\b_y(x)\bigr) \ \ (x\in B').
$$
Moreover, for $y\in Y_0$ (when the sections $a_y$ and $b_y$ agree over $\tilde C$
and hence define a section over a \nbd\ of $A\cup B$) we have by construction
$\psi_y(x,0)=0$ for all $x\in \tilde C$. Therefore $\a_y=0$ and $\b_y=0$ for
such $y$ (prop.\ 5.2), and hence the above homotopies are fixed for $y\in Y_0$
as required. Everything else is clear from the earlier arguments.
\endpr

%
%
\proclaim 5.6 Corollary: Let $X$ be a Stein manifold and
$A\subset A'$ a non-critical \spsc\ extension in $X$ (def.\ 2.1).
Let $F$ be a complex manifold with a spray
and $Z=X\times F\to X$ the associated trivial bundle with fiber $F$. Let
$d$ be a metric on $Z$. Given a homotopy $a_t$ ($0\le t\le 1$) of holomorphic
sections of $Z$ in a \nbd\ of $A$ such that $a_0$ and $a_1$ are holomorphic
in a \nbd\ of $A'$, there is for each $\e>0$ a homotopy $\tilde a_t$ ($0\le t\le 1$)
of holomorphic sections in a \nbd\ of $A'$ such that $\tilde a_0=a_0$,
$\tilde a_1=a_1$, and
$$
    d\bigl( \tilde a_t(x), a_t(x) \bigr) <
    \e \quad (x\in A,\ 0\le t\le 1).
$$

\ni\it Proof. \rm
By theorem 2.5 there is a finite sequence
$$ A=A_0\subset A_1\subset A_2\subset\ldots \subset A_k  = A' $$
such that for each $j=0,1,\ldots, k-1$ we have $A_{j+1}=A_j\cup B_j$,
where $(A_j,B_j)$ is a convex bump. Suppose inductively that
we have already approximated the initial homotopy $a_t=a^0_t$ on $A$
by a homotopy $a^j_t$ ($0\le t\le 1$) consisting of sections that are holomorphic in
a \nbd\ of $A_j$ (with $a^j_0=a_0$ and $a^j_1=a_1$). Since $B_j$ is a convex bump
on $A_j$, we can apply theorem 4.5 to approximate the homotopy
$a^j_t$ in a \nbd\ of $C_j=A_j\cap B_j$ by a homotopy $b_t$ consisting
of holomorphic sections in a \nbd\ of $B_j$ and such that
$b_0=a_0$, $b_1=a_1$. If the approximation is sufficiently close,
we can apply theorem 5.5 to glue the families $a^j_t$ and $b_t$
over $C_j$ into a single family $a^{j+1}_t$ ($0\le t\le 1$)
consisting of sections which are holomorphic in a \nbd\ of
$A_{j+1}=A_j\cup B_j$, which equal $a_0$ resp.\ $a_1$ when
$t=0$ resp.\ $t=1$, and which approximate the sections $a^j_t$
on $A_j$. This completes the induction step. In a finite number
of steps we reach a desired homotopy $\tilde a_t=a^k_t$
satisfying corollary 5.6.
\endpr

\demo Remark: The proof of corollary 5.6 remains valid if the bundle $Z$ is trivial
over a \nbd\ of each bump $B_j$ but not necessarily over $A'$. Thus the corollary
provides an extension of the h-principle from $A$ to any non-critical
\spsc\ extension of $A$.

%
%
%
%
\beginsection 6. Extending holomorphic sections across \pb s.

In this section we prove theorem 2.6, thereby concluding the proof of
theorems 1.3 and 1.4. We present the proof in the case without
parameters, but all arguments go through in the general parametric
case by using suitable versions of the approximation and gluing
theorems which were proved in sections 4 and 5.
There are three main steps:

\medskip
\item{--} approximate the given section $a$ in a \nbd\ of $C=A\cap B$
by a section $b$ which is holomorphic in a \nbd\ of $B$;

\item{--} glue the sections $a$ and $b$ by theorem 5.1 (or theorem 5.5)
to obtain a holomorphic section $\tilde a$ in a \nbd\ of $A\cup B$;

\item{--} show that the new section $\tilde a$ is homotopic to $a$ and satisfies
all required properties.

\medskip
We present all arguments in the general case when $(A,B)$ is a \psc\ bump
and the set $C=A\cap B$ is not necessarily contractible. In the special case of
convex bumps the arguments involving homotopies can be substantially simplied
by using contractibility of $C$ and the fact that $A$ is a strong deformation
retraction of $A\cup B$.

\proclaim 6.1 Proposition: There exists a homotopy $b_t$  ($0\le t\le 1$)
of continuous sections in a \nbd\ of $B$ such that $b_0$ is holomorphic
in a \nbd\ of $B$, each $b_t$ is holomorphic in a \nbd\ of $C$, and
$b_1=a$.

\demo Proof: Since $B$ is star-shaped, it has an open contractible \nbd\ $\tilde B$.
By hypothesis we may choose $\tilde B$ so small that $Z|\tilde B$ is a trivial bundle.
Hence there is a holomorphic section $b_0$ of $Z$ over $\tilde B$ (we may simply
take a constant section in a given trivialization) and a homotopy
$b'_t \colon \tilde B \to Z$ ($0\le t\le 1$) of continuous sections connecting
$b'_0=b_0$ and $b'_1=a$ (such a homotopy exists since $\tilde B$ is contractible
and the bundle is trivial over $\tilde B$).

Recall (part (v) in def.\ 2.2) that we have a \spsh\ function $\tau \ge 0$
in a \nbd\ of $C=A\cap B$ such that $S=\{\tau=0\}$ is a \tr\ sphere
contained in an affine plane $T$ (in some holomorphic coordinates
in a \nbd\ of $B$), $C=\{\tau\le 1\}$, and $\tau$ has no critical
points on $C\bs S$. Our first goal is to modify the
homotopy $b'_t$ so as to make it holomorphic in a \nbd\ of $S$.
Assume (as we may) that $S$ is \ra.  Let $T'$ be the affine \tr\ subspace
of maximal dimension containing $T$ (and hence $S$). We can approximate $b'_t|S$
uniformly on $S$ by a homotopy consisting of \ra\ sections over a \nbd\ of $B$
without changing the sections $b'_0=b_0$ and $b'_1=a$. We still denote this
\ra\ homotopy by $b'_t$. For each $t$ the section $b'_t|T'$ extends (by complexification)
to a unique holomorphic section $\tilde b_t$ in a \nbd\ of $T'\cap C$
(independent of $t$). Hence $\tilde b_t$ is a holomorphic homotopy in a \nbd\ of
$S_0= \{\tau \le  c_0\}$ for some sufficiently small $c_0>0$.   Of course this process
does not affect the sections that were already holomorphic, so we have
$\tilde b_0=b_0$ and $\tilde b_1=a$ on $S_0$.

Since $C$ is a non-critical \spsc\ extension of $S_0$, corollary 5.6 implies that
$\tilde b_t$ can be approximated uniformly on $S_0$ by another homotopy $b_t$
($0\le t\le 1$) which is holomorphic in a \nbd\ of $C$ and connects $b_0$ and $a$.

To complete the proof we must show that the homotopy $b_t$ (which has so far been defined
and holomorphic in a \nbd\ of $C$) extends to a continuous homotopy from
$b_0$ to $a$ in a \nbd\ of $B$. To do this we first reparametrize both homotopies
$b_t$ and $b'_t$ so that for some small $\d>0$ we have
$$ b_t=b'_t=b_0 \quad {\rm for\ } 0\le t\le \d,\qquad
   b_t=b'_t=a \quad {\rm for\ }   1-\d \le t\le 1 .         \eqno(6.1) $$

\proclaim 6.2 Lemma: If $b_t|S$ is sufficiently uniformly close to $b'_t|S$ on
$S=\{\tau=0\} \subset C$ for each $t\in [0,1]$, and if (6.1) holds for some
$\d>0$, then there exist a \nbd\ $\tilde C$ of $C$ and a two-parameter homotopy
$\tilde b_{t,s}$ ($0\le t,s\le 1$) of continuous sections of $Z$ over $\tilde C$
satisfying
$$ \eqalign{ & \tilde b_{t,0}=b_t,\quad \tilde b_{t,1}=b'_t\ \
                                  {\rm for}\ \ 0\le t\le 1; \cr
    & \tilde b_{t,s}=b_t=b'_t\ \ {\rm for}\ \  t\in \{0,1\}
    \ \ {\rm and}\ \ 0\le s\le 1. \cr}
$$

\demo Proof:  If $b_t|S$ is sufficiently uniformly close to $b'_t|S$ for each
$t \in [0,1]$ (which we may assume to be the case), we can use the spray as in
lemma 3.2 to obtain a two parameter homotopy of sections
$b_{t,s} \colon S\to Z$ ($0\le t,s\le 1$) satisfying
$$ \eqalign{ & b_{t,0}=b_t|S,\quad b_{t,1}=b'_t|S\ \ {\rm for}\ \ 0\le t\le 1; \cr
    & b_{t,s}=b_t=b'_t\ \ {\rm for}\ \  t\in [0,\d] \cup [1-\d,1]
    \ \ {\rm and}\ \ 0\le s\le 1. \cr}
$$
Furthermore, $S$ is a strong deformation retraction of a \nbd\ $\tilde C$ of $C$,
i.e., there is a smooth family of maps
$\Theta_t\colon \tilde C\to \tilde C$ ($0\le t\le 1$) such that $\Theta_0$ is the
identity on $\tilde C$, each $\Theta_t$ is the identity on $S$, and
$\Theta_1(\tilde C) =S$. Choose a smooth function $u\colon [0,1]\to [0,1]$
which is zero in a \nbd\ of $t=0$ and $t=1$ and is identically one on
$[\d,1-\d]$. We identify sections of $Z$ over $\tilde C$ (or over $\tilde B$)
by mappings into the fiber $F$, using the triviality of $Z|\tilde B$.
For each $x\in \tilde C$ and $0\le t\le 1$ we define
$$  \tilde b_{t,s}(x) = \cases{
       b_{t}\bigl( \Theta_{3su(t)} (x)\bigr);   & if $0\le s\le 1/3$,\cr
       b_{t,3s-1}\bigl( \Theta_{u(t)}(x)\bigr); & if $1/3\le s\le 2/3$, \cr
       b'_t\bigl( \Theta_{3(1-s)u(t)}(x)\bigr); & if $2/3\le s\le 1$.\cr}
$$
Clearly these maps are defined for all $x\in \tilde C$. We have
$\tilde b_{t,0}=b_t\circ\Theta_0=b_t$ and $\tilde b_{t,1}=b'_t\circ\Theta_0=b'_t$.
Moreover, when $t$ is near $0$ or near $1$, or when $x\in S$, we have
$\tilde b_{t,s}(x)= b_{t,s}(x)$ for all $s\in [0,1]$. Thus $\tilde b_{t,s}$
satisfies lemma 6.2.
\endpr

Choose a smooth function $v\colon X\to  [0,1]$ such that $v=1$
in a \nbd\ $C' \subset \tilde C$ of $C$ and $\supp v\ss \tilde C$.
Let $\tilde b_{t,s}$ be as in lemma 6.2. Consider the homotopy
$$ \tilde b_{t,1-v(x)}(x) \quad (x\in \tilde B,\ 0\le t\le 1). $$
By the choice of $v$ this equals $\tilde b_{t,0}(x)= b_t(x)$ for $x\in C'$
(so it is holomorphic there), it is defined for all $x\in \tilde B$,
and it equals $b'_t(x)$ for $x$ outside $\tilde C$. This homotopy, which we
again denote by $b_t$, satisfies proposition 6.1.
\endpr

In what follows we shall shrink the neighborhoods $\tilde B$ and $\tilde C$
several times without mentioning this again, and without changing the notation.
We are now in position to apply the h-Runge approximation (theorem 4.1 or 4.2)
to approximate the homotopy $b_t$ from proposition 6.1 uniformly in
$\tilde C$ by a holomorphic homotopy $\tilde b_t\colon \tilde B\to Z$ such that
$\tilde b_0=b_0$. The section $b=\tilde b_1$ then approximates $a$ in $\tilde C$
as well as desired. Hence by theorem 5.1 (or 5.5) we can glue $a$ and $b$ into
a single section $\tilde a$ which is holomorphic in a \nbd\ $A'$ of $A\cup B$
and which approximates $a$ in a \nbd\ of $A$.

It remains to show that there is a homotopy of sections $a_t$ ($0\le t\le 1$)
in a \nbd\ of $A\cup B$, connecting $a_0=a$ and $a_1=\tilde a$, such that
each $a_t$ is holomorphic in a \nbd\ of $A$ and approximates $a$ there
(so $a_t$ will satisfy theorem 2.6). In a \nbd\ of $A$ such a homotopy
is provided by theorem 5.1. Our goal is to extend this
homotopy to a \nbd\ of $A\cup B$ by modifying it outside some \nbd\ of $A$.

Recall that over $\tilde B$ we have the following homotopies:

\item{(i)}  the homotopy $b_t$ from $b_0$ to $b_1=a$,
given by proposition 6.1;

\item{(ii)} the homotopy $\tilde b_t$ from $\tilde b_0=b_0$
to $\tilde b_1=b$, obtained by approximating $b_t$ in $\tilde C$;

\item{(iii)} the homotopy $b'_t$ from $b'_0=b(=\tilde b_1)$ to
$b'_1=\tilde a$, given by theorem 5.1.

Note that all these homotopies are holomorphic in $\tilde C$ and $b'_t$
approximates $b$ (and hence $a$) there.
If we combine these three homotopies in the correct order (first follow
$b_{1-t}$ from $a$ to $b_0$, then follow $\tilde b_t$ from $b_0$ to $\tilde b_1=b$,
and finally follow $b'_t$ from $b$ to $\tilde a$), we get a homotopy from
$a$ to $\tilde a$ over $\tilde B$ which is holomorphic over $\tilde C$. However,
we must show that this homotopy over $\tilde B$ can be glued with the homotopy
$a_t$ into a single homotopy from $a$ to $\tilde a$ in a \nbd\ of $A\cup B$.

In order to do this we will first join the above homotopies (i)--(iii) over $\tilde B$
into a new homotopy $h_t$ from $a$ to $\tilde a$ which in addition will
approximates $a$ in $\tilde C$. For convenience we shall define $h_t$ initially
on the $t$-interval $[0,4]$ and subsequently rescale the parameter to $[0,1]$.
According to the remark following theorem 4.1 we may assume that there
is a two-parameter homotopy $g_{t,s} \colon \tilde B\to Z$\ ($0\le t,s\le 1$)
which is holomorphic in $\tilde C$ and satisfies
$$ g_{t,0}=b_t,\quad g_{t,1}=\tilde b_t,\quad g_{0,s}=b_0, $$
and such that $g_{1,s}$ approximates $b_1=a$ in $\tilde C$ for each $s\in [0,1]$.
Choose a smooth function $\chi$ on $X$ with values in $[0,1]$ such that
$\chi=1$ in a \nbd\ of $A$ and $\chi=0$ in a \nbd\ of $B\bs \tilde C$.
Such $\chi$ exists since the sets $A\bs B$ and $B\bs A$ are separated.
We now define for each $x\in \tilde B$
$$     h_t(x)= \cases{
            b_{1-t+t\chi(x)}(x), &if $0\le t\le 1$; \cr
            g_{\chi(x),t-1}(x), &if $1\le t\le 2$;\cr
            \tilde b_{\chi(x)+(1-\chi(x))(t-2)} (x), &if  $2\le t\le 3$;\cr
            b'_{t-3}(x), &if $3\le t\le 4$. \cr}
$$
The reader may verify that this is indeed a homotopy from $h_0=a$ to
$h_4=\tilde a$. For $x$ in a \nbd\ of $C$ we have $\chi(x)=1$ and hence
$$      h_t(x)= \cases{
            b_1(x)=a(x),   &if $0\le t\le 1$; \cr
            g_{1,t-1}(x),  &if $1\le t\le 2$;\cr
            \tilde b_1(x)=b(x), &if  $2\le t\le 3$;\cr
            b'_{t-3}(x),   &if $3\le t\le 4$. \cr}
$$
Hence $h_t$ is holomorphic near $C$ and approximates $a$ there (since all
homotopies in question are close to $a$ on $\tilde C$). Hence we may assume that
on $\tilde C$ both homotopies $a_t$ and $h_t$ (which we rescale to the $t$-interval
$[0,1]$) approximate $a$ so well that their images (over $\tilde C$) belong to
a tubular \nbd\ of $a(\tilde C) \subset Z$ in which we can apply lemma 3.2.
This means that we can view these sections on $\tilde C$ as sections of a certain
holomorphic vector bundle over $\tilde C$. This allows us to find a two parameter
homotopy $k_{t,s}$ joining $a_t$ and $h_t$ over $\tilde C$ (we can simply use
the convex combinations of the two sections in the given vector bundle.) Finally
we patch $a_t$ and $h_t$ using $k_{t,s}$ into the homotopy
$$ \tilde a_t(x) = k_{t,1-\chi(x)} (x)                                  $$
where $\chi$ is a smooth function chosen as above. For $x$ near $A$ we
have $\chi(x)=1$ and hence $\tilde a_t(x)=k_{t,0}(x)=a_t(x)$, while for $x$
near $B\bs \tilde A$ we have $\chi(x)=0$ ad hence
$\tilde a_t(x)=k_{t,1}(x)=h_t(x)$.  We denote this new homotopy again $a_t$.

Finally we choose a smooth function $\eta$ on $X$ with values in $[0,1]$
such that $\eta=1$ near $A\cup B$ and $\supp \eta \subset A'$, where $A'$
is the \nbd\ of $A\cup B$ on which the homotopy $a_t$ has been defined.
The homotopy $a_{t\eta(x)}(x)$ ($0\le t\le 1$) is now defined
for all $x\in X$, it equals $a_0(x)=a(x)$ for $x\in X\bs A'$ and $t\in [0,1]$, and
it equals $a_t(x)$ for $x\in A\cup B$. This completes the proof of theorem 2.6.

%
%
%
%
\medskip\ni\bf References. \rm \medskip

\ii{[AFR]} P.\ Ahern, M.\ Flores, J.-P.\ Rosay:
On $\R^+$ and $\C$ complete holomorphic vector fields.
Proc.\ Amer.\ Math.\ Soc., to appear.

\ii{[Car]} H.\ Cartan: Espaces fibr\'es analytiques.
Symposium Internat.\ de topologia algebraica, Mexico, 97--121 (1958).
(Also in Oeuvres, vol.\ 2, Springer, New York, 1979.)

\ii{[BFo]} G.\ Buzzard, J.\ E.\ Forn\ae ss:
An embedding of $\C$ in $\C^2$ with hyperbolic complement.
Math.\ Ann.\ {\bf 306}, 539--546 (1996).

\ii{[Chi]} E.\ M.\ Chirka: Complex analytic sets.
Mathematics and its Applications (Soviet Series), {\bf 46},
Kluwer, Dordrecht, 1989.

\ii{[Dem]} J.-P.\ Demailly:
Un exemple de fibr\'e holomorphe non de Stein \`a fibre $\bf C^2$ ayant pour
base le disque ou le plan.
Invent.\ Math.\ {\bf 48}, 293--302 (1978).

\ii{[Eli]} Y.\ Eliashberg:
Topological characterization of Stein manifolds of dimension $>2$.
Internat.\ J.\ Math.\ {\bf 1}, 29--46 (1990).

\ii{[EGr]} Y.\ Eliashberg, M.\ Gromov: Embeddings of Stein manifolds.
Ann.\ Math.\ {\bf 136}, 123--135 (1992).

\ii{[For]} O.\ Forster:
Topologische Methoden in der Theorie der Steinscher R\"aume.
(Iternat.\ Congress in Math., Nice, 1970, pp.\ 613--618)
Gauthier-Villars, Paris, 1971.

\ii{[FR1]} O.\ Forster and K.\ J.\ Ramspott:
Okasche Paare von Garben nicht-abelscher Gruppen.
Invent.\ Math.\ {\bf 1}, 260--286 (1966).

\ii{[FR2]} O.\ Forster and K.\ J.\ Ramspott:
Analytische Modulgarben und Endromisb\"undel.
Invent.\ Math.\ {\bf 2}, 145--170 (1966).

\ii{[Fo1]} F.\ Forstneri\v c:
Actions of $(\R,+)$ and $(\C,+)$ on complex manifolds.
Math.\ Z.\ {\bf 223}, 123--153 (1996).

\ii{[Fo2]} F.\ Forstneri\v c:
Interpolation by holomorphic automorphisms and embeddings in $\C^n$.
J.\ Geom.\ Anal.\ {\bf 9}, no.1, (1999) 93-118.

\ii{[FGR]} F.\ Forstneri\v c, J.\ Globevnik, J.-P.\ Rosay:
Non straightenable complex lines in $\C^2$.
Arkiv Mat.\ {\bf 34}, 97--101 (1996).

\ii{[FP]} F.\ Forstneri\v c and J.\ Prezelj:
Oka's principle for holomorphic submersions with sprays.
Preprint, 1999.

\ii{[Gr1]} H.\ Grauert:
Charakterisierung der holomorph vollst\"andigen R\"aume.
Math.\ Ann.\ {\bf 129}, 233--259 (1955).

\ii{[Gr2]} H.\ Grauert:
Approximationss\"atze f\"ur holomorphe Funktionen mit Werten in
komplexen R\"aumen.
Math.\ Ann.\ {\bf 133}, 139--159 (1957).

\ii{[Gr3]} H.\ Grauert:
Holomorphe Funktionen mit Werten in komplexen Lieschen Gruppen.
Math.\ Ann.\ {\bf 133}, 450--472 (1957).

\ii{[Gr4]} H.\ Grauert: Analytische Faserungen über
holomorph-vollständigen Räumen.
Math.\ Ann.\ {\bf 135}, 263--273 (1958).

\ii{[GRe]} H.\ Grauert, R.\ Remmert: Theory of Stein Spaces.
Grundl.\ Math.\ Wiss.\ {\bf 227}, Springer, New York, 1977.

\ii{[Gro]} M.\ Gromov:
Oka's principle for holomorphic sections of elliptic bundles.
J.\ Amer.\ Math.\ Soc.\ {\bf 2}, 851-897 (1989).

\ii{[GRo]} C.\ Gunning, H.\ Rossi:
Analytic functions of several complex variables.
Prentice--Hall, Englewood Cliffs, 1965.

\ii{[HKu]} P.\ Heinzner and F.\ Kutzschebauch:
An equivariant version of Grauert's Oka principle.
Invent.\ Math.\ {\bf 119}, 317--346 (1995).

\ii{[HL1]} G.\ M.\ Henkin, J.\ Leiterer:
Theory of functions on complex manifolds.
Akademie-Verlag, Berlin, 1984.

\ii{[HL2]}  G.\ Henkin, J.\ Leiterer:
Proof of Oka-Grauert principle without the induction over basis dimension.
Preprint, Karl Weierstrass Institut f\"ur Mathematik,
Berlin, 1986.

\ii{[HL3]} G.\ Henkin, J.\ Leiterer:
The Oka-Grauert principle without induction over the basis dimension.
Math.\ Ann.\ {\bf 311}, 71--93 (1998).

\ii{[H\"o1]} L.\ H\"ormander:
$L\sp{2}$ estimates and existence theorems for the $\dibar$ operator.
Acta Math.\ {\bf 113}, 89--152 (1965).

\ii{[H\"o2]} L.\ H\"ormander:
An Introduction to Complex Analysis in Several Variables, 3rd ed.
North Holland, Amsterdam, 1990.

\ii{[Lei]} J.\ Leiterer:
Holomorphic Vector Bundles and the Oka-Grauert Priciple.
Encyclopedia of Mathematical Sciences, vol. 10, 63--103;
Several Complex Variables IV, Springer,  1989.

\ii{[Ram]} K.\ J.\ Ramspott:
Stetige und holomorphe Schnitte in Bundeln mit homogener Faser.
Math.\ Z.\ {\bf 89}, 234--246 (1965).

\ii{[Ros]} J.-P.\ Rosay: A counterexample related to Hartogs
phenomenon (a question by E. Chirka).
Michigan Math.\ J.\ {\bf 45}, 529--535 (1998).

\ii{[RRu]} J.-P.\ Rosay, W.\ Rudin:
Holomorphic maps from $\C^n$ to $\C^n$.
Trans.\ Amer.\ Math.\ Soc.\ {\bf 310}, 47--86 (1988)

\ii{[Shd]} M.\ Schneider: Tubenumgebungen Steinscher R\"aume.
Manuscripta Math.\ {\bf 18}, 391--397 (1976).

\ii{[Sch]} J.\ Sch\"urmann:
Embeddings of Stein spaces into affine spaces of minimal dimension.
Math.\ Ann.\ {\bf 307}, 381--399 (1997).

\ii{[Siu]} Y.T.\ Siu: Every Stein subvariety admits a Stein neighborhood.
Invent.\ Math.\ {\bf 38}, 89--100 (1976).

\ii{[Ste]} K.\ Stein: Analytische Funktionen mehrerer komplexer
Veränderlichen zu vorge\-ge\-benen Periodizitätsmoduln und das zweite
Cousinsche Problem.
Math.\ Ann.\ {\bf 123}, 201--222 (1951).

%
%
%
%
\bigskip\medskip
\settabs 5\columns
\+\ \ Franc Forstneri\v c            &&& Jasna Prezelj\cr
\+\ \ Department of Mathematics      &&& Faculty of Mechanical Engineering \cr
\+\ \ University of Wisconsin        &&& University of Ljubljana \cr
\+\ \ Madison, WI 53706, USA         &&& A\v sker\v ceva 6\cr
\+                               &&& 1000--Ljubljana, Slovenia \cr
\+\ \ \it Current address:\rm \cr
\+\ \ IMFM, University of Ljubljana  \cr
\+\ \ Jadranska 19  \cr
\+\ \ 1000--Ljubljana, Slovenia\cr

\bye